\documentclass{amsart}
\usepackage{latexsym, amssymb, amsmath, amsthm, rotating, amscd}
\usepackage{graphics}
\newtheorem{theorem}{Theorem}[section]
\newtheorem{definition}[theorem]{Definition}

\newtheorem{lemma}[theorem]{Lemma}
\newtheorem{proposition}[theorem]{Proposition}

\newtheorem{remark}[theorem]{Remark}
\newtheorem{example}[theorem]{Example}
\numberwithin{equation}{section}
\newenvironment{ack}[1][Acknowledgements]{\begin{trivlist}
\item[\hskip \labelsep {\bfseries #1}]}{\end{trivlist}}

\def\ep{\varepsilon}
\def\Z{\mathbb{Z}}
\def\N{\mathbb{N}}
\def\C{\mathbb{C}}
\def\R{\mathbb{R}}
\def\co{\begin{sideways}\begin{sideways}$Y$\end{sideways}\end{sideways}}

\begin{document}

\title[Vertex operator coalgebras]{The notion of vertex operator
coalgebra and a geometric interpretation}
\author{Keith Hubbard}
\date{\today}

\begin{abstract}
The notion of vertex operator coalgebra is presented and motivated
via the geometry of conformal field theory. Specifically, we
describe the category of geometric vertex operator coalgebras,
whose objects have comultiplicative structures meromorphically
induced by conformal equivalence classes of worldsheets. We then
show this category is isomorphism to the category of vertex
operator coalgebras, which is defined in the language of formal
algebra.  The latter has several characteristics which give it the
flavor of a coalgebra with respect to the structure of a vertex
operator algebra and several characteristics that distinguish it
from a standard dual -- both will be mentioned.
\end{abstract}

\maketitle

\section{Introduction}

In this paper, we define the notion of a vertex operator coalgebra
as motivated through the geometry of conformal field theory.
Specifically, after defining the category whose objects are vertex
operator coalgebras (VOCs) we show that this category is
isomorphic to the category whose objects are geometrically defined
structures called geometric vertex operator coalgebras (GVOCs).
The results presented in this paper reflect a portion of the
results proven in the author's Ph.D. dissertation \cite{K}. The
definition of a geometric vertex operator coalgebra is motivated
in a similar way to that of a geometric vertex operator algebra
(\cite{H2}, \cite{H}) in that both arise from the geometry
underlying conformal field theory, specifically closed string
theory, in physics.  Conceptually, string theory seeks to unify
all the fundamental forces in nature into a single theory by
modelling particles as one-dimensional objects (as opposed to
point-particles) moving through space-time.  These objects, called
strings, sweep out two-dimensional manifolds, called worldsheets,
in space-time that model the interactions of particles.

Closed string theory focuses specifically on strings that begin
and end at the same point, that is, on loops. Their study, along
with the study of the algebraic structure they induce began in
physics (cf. \cite{BPZ}, \cite{FS}, \cite{S} and \cite{V}).
Physicists introduced the notion of vertex operators in order to
write down the expectation values of certain particle
interactions.  These vertex operators were independently
discovered in mathematics in the study of representations of
affine Lie algebras (\cite{LW}).  Shortly afterwards, Borcherds,
attempting to describe the algebra arising from particle states
and their vertex operators in a concise mathematical structure,
introduced the notion of a vertex algebra (\cite{Bo}).  This
notion was then specialized in \cite{FLM} to the notion of a
vertex operator algebra (VOA).  Although vertex operator algebras
had origins in the geometry of string theory, a rigorous
correspondence between the algebraic formalism and the geometric
motivation was not completed until Huang introduced the notion of
a geometric vertex operator algebra (GVOA) in \cite{H2} and
\cite{H}.  Motivated by the physicists' study of the geometric
interaction (specifically chiral interaction) of $n \in \N$ closed
strings combining to form one closed string, Huang defined
geometric vertex operator algebras in the context of a
mathematical description of genus-zero Riemann surfaces with tubes
up to conformal equivalence, along with a sewing operation on
these objects.  He also proved that the category of GVOAs was
isomorphic to the category of VOAs, finalizing the rigorous
correspondence between VOAs and the geometry of genus-zero
worldsheets. The interactions corresponding to the worldsheets
that Huang focused on are important special cases of interactions
that may occur in a conformal field theory.

Another way that Huang's definition of GVOAs proved fruitful, was
that it showed a VOA could be described in the language of
operads, a language first developed by May in \cite{M} to discuss
iterated loop spaces.  Specifically, Huang's GVOA corresponds to
describing the partial operad of equivalence classes of genus-zero
compact Riemann surfaces with $n$ incoming punctures, one outgoing
puncture and local coordinates vanishing at each puncture, then
taking the meromorphic associative algebra associated to
$\C$-extensions of that partial operad (\cite{HL}). Using this
interpretation, attaching or ``sewing" two worldsheets together
corresponds to composing operators in an operad defined in the
category of complex vector spaces.

Prior to VOAs being solidly tied to their foundations in the
geometry of conformal field theory, significant applications for
VOAs had been emerging outside of conformal field theory. VOAs are
now known to have deep connections to modular functions, have
contributed to the representation theory of the Monster finite
simple group and are fundamental to the study of representations
of infinite-dimensional Lie algebras. (See the introduction of
\cite{FLM} for a more thorough exposition on connections and
applications.)

Geometric VOAs are quite effective in modelling the interactions
for which they were designed, specifically the combining of $n$
incoming strings into one outgoing string.  There are, however,
still very significant questions about what occurs in other
settings.  Geometric vertex operator coalgebras are designed to
model genus-zero worldsheets corresponding to one incoming closed
string splitting into $n$ outgoing closed strings. Similar to the
GVOA case, this amounts to describing the partial operad of
equivalence classes of genus-zero compact Riemann surfaces with
one incoming and $n$ outgoing punctures, then investigating the
algebraic structure it induces on complex vector spaces, i.e.
determining the meaning of a coalgebra over this partial operad.
The interactions modelled by GVOCs are of interest in their own
right but also open the possibility of describing even more
particle interactions when combined with existing GVOA theory. For
instance, sewing the one incoming puncture of a worldsheet for a
GVOC with the one outgoing puncture of a worldsheet for a GVOA,
one obtains a worldsheet with any number of incoming and any
number of outgoing punctures.  Attaching the incoming punctures of
a worldsheet associated to a GVOA with the outgoing punctures of a
worldsheet associated to a GVOC forms a higher-genus worldsheet.
Exploring the algebraic compositions of the two operators
corresponding to these two worldsheets may shed light on
understanding higher-genus conformal field theory, which has been
an area of marked difficulty in the past.  It is also significant
to note that any full conformal field theory must include VOCs as
well as VOAs.

The main thrust of this paper, however, is to understand VOCs from
an algebraic and geometric standpoint.  To that end, the notion of
a vertex operator coalgebra is defined (purely in terms of vector
spaces, formal variables and formal operators) and is shown to be
categorically isomorphic to the notion of GVOC.  Generating a
formal algebraic description in the VOA case has facilitated a
more concrete understanding of the structure as well as the
construction of examples. Similarly the algebraic description of
VOCs allows us to describe concrete examples (\cite{K2}) and opens
the door to a deeper understanding of a conformal field theory's
algebraic structure.  One might reasonably wonder whether the
definition of VOC might be arrived at through purely algebraic
means.  But ``standard" dualizing would not produce bounding of
formal variables in opposite directions that arise in VOCs -
specifically, truncation produces only finitely many positive
powers of the given formal variable in a VOC, while the counit
property produces only non-negative powers (in VOAs everything is
truncated from below). Additionally, results of weak commutativity
coming from the VOC Jacobi identity reveal an important
distinction from that of the VOA Jacobi identity (cf. \cite{K3}).

Since VOAs are of interest outside of string theory, it is natural
to hope that VOCs, once well understood, might introduce new
connections to broader mathematics as well.

In Section 2 we will recall Huang's conventions (\cite{H2},
\cite{H}) for spheres with multiple outgoing tubes and make the
obvious generalizations of notation and theory to spheres with
multiple incoming tubes.  Section 3 will cover the necessary
formal calculus to describe the context in which our
comultiplication operators will be defined.  In Section 4 we will
define the notion of geometric vertex operator coalgebra as well
as the notion of (algebraic) vertex operator coalgebra and discuss
why finding a coalgebra structure is more subtle than the process
of reversing the orientation of spheres with tubes (see Remark
\ref{R:subtle}). There will also be a number of properties of VOCs
presented in Section 4. Section 5 will provide a constructive
isomorphism from the category of GVOCs to the category of VOCs as
well as its inverse.

\begin{ack}
This paper reflects a portion of the author's Ph.D. dissertation.
The author would like to thank his advisor, Katrina Barron, for
her guidance and encouragement throughout the research, writing
and revising phases of this paper and his dissertation.  In
addition, he would like to express thanks to Stephan Stolz for
numerous helpful conversations.  The author also gratefully
acknowledges the financial support of the Arthur J. Schmidt
Foundation.
\end{ack}

\section{The geometry of spheres with tubes} \label{S:geometry}
We begin by investigating the geometry of spheres with tubes.
Since our main motivation in this paper is studying closed strings
propagating through space-time, it is essential that we precisely
define the types of structures that model this behavior.  The
description in this section will closely follow \cite{H} since
that work contains the description of the geometry underlying
GVOAs (motivated by string theory) and also because we would like
the geometry of GVOCs to be compatible with that structure.

\subsection{Defining spheres with tubes} \label{S:define}
Geometric vertex operator algebras and coalgebras are defined
using the moduli space (under conformal equivalence) of genus-zero
compact Riemann surfaces with punctures and local coordinates
vanishing at each puncture.  Vafa first observed that having a
puncture on a Riemann surface together with local coordinates that
vanish at that puncture is conformally equivalent to having a
half-infinite tube attached to the Riemann surface (cf. \cite{V},
\cite{H}).  We will begin by defining exactly what we mean by a
genus-zero Riemann surface with punctures and local coordinates
vanishing at the punctures, or equivalently, a genus-zero
worldsheet.

Consider a compact genus-zero Riemann surface, by which we mean a
compact connected genus-zero one-dimensional complex manifold.  In
\cite{H} and in this paper these surfaces are simply referred to
as ``spheres".  The reason for the use of the word spheres is that
any compact Riemann surface of genus-zero is complex analytically
isomorphic to the standard Riemann sphere, i.e. $\C \cup \{ \infty
\}$ with the standard complex structure  (cf. \cite{A}, \cite{F}).
Since we will eventually be concerned only with conformal
equivalence classes of Riemann surfaces, in Section
\ref{S:modulispace} we will pick canonical representatives of the
equivalences classes under conformal equivalence and these
canonical representatives will be Riemann spheres.

By the term \emph{oriented puncture}
we mean the selection of a point of a sphere together with an element of
\{+, -\}.
Given an oriented puncture $p$, a \emph{local analytic coordinate chart vanishing at $p$}
is a pair $( U, \Phi )$ where
$U$ is an open neighborhood of $p$ in the sphere, called the \emph{local coordinate neighborhood}, and
$\Phi : U \to \C$ is an injective
analytic map, such that $\Phi (p) = 0$, called the \emph{local coordinate map}.  The term
\emph{tube centered at $p$} is used interchangeably with the term local analytic coordinate
chart vanishing at $p$ and may by thought of as the path of a single closed string propagating through
space-time.  Tubes centered at negatively oriented punctures represent half-infinite tubes swept out
by outgoing strings, while tubes centered at positively oriented punctures represent half-infinite
tubes swept out by incoming strings.

Spheres with $m$ tubes centered at negatively oriented punctures and $n$ tubes centered at
positively oriented punctures are said to be of \emph{type} $(m,n)$.  All oriented punctures are
required to be distinct and have an ordering with negatively oriented punctures coming first.
In \cite{H}, the focus is primarily on spheres of type $(1,m)$, with $m \in \N$.  In this work, we will
discuss spheres of type $(1,m)$, for $m \in \N$,  called \emph{spheres with incoming tubes}, and
of type $(m,1)$ for $m \in \N$,  called \emph{spheres with outgoing tubes}, but will eventually focus
on the latter.

For the moment, however, we will discuss all genus-zero worldsheets.  Using the above framework,
we can denote a sphere (or compact Riemann surface of genus-zero) with
$m$ outgoing punctures and $n$ incoming punctures by

\begin{equation*}
(S; p_{-m}, \ldots, p_{-1},p_{1}, \ldots, p_{n}; (U_{-m},\Phi_{-m}), \ldots, (U_{-1},\Phi_{-1}),
(U_{1},\Phi_{1}), \ldots, (U_{n},\Phi_{n}))
\end{equation*}

\noindent
where $m,n \in \N$, $S$ is a compact genus-zero Riemann surface, $p_i$ is a point in $S$ with the
sign of the index corresponding to the orientation of the puncture,
and $(U_i,\Phi_i)$ is a local analytic coordinate chart vanishing at $p_i$ for $i=-m, \ldots -1,1,
\ldots,n$.  The terms \emph{sphere with tubes}, \emph{genus-zero worldsheet}, and
\emph{compact Riemann surface of genus-zero with ordered punctures and local coordinates vanishing at
the punctures} are used interchangeably
in the literature to refer to this structure.

Let
\begin{multline*}
\Sigma_1=(S_1; p_{-m}, \ldots, p_{-1},p_{1}, \ldots, p_{n}; \\
(U_{-m},\Phi_{-m}), \ldots, (U_{-1},\Phi_{-1}),
(U_{1},\Phi_{1}), \ldots, (U_{n},\Phi_{n}))
\end{multline*}
\noindent
be a sphere of type $(m,n)$ and
\begin{multline*}
\Sigma_2=(S_2; q_{-k}, \ldots, q_{-1},q_{1}, \ldots, q_{\ell}; \\
(V_{-k},\Psi_{-k}), \ldots, (V_{-1},\Psi_{-1}),
(V_{1},\Psi_{1}), \ldots, (V_{\ell},\Psi_{\ell}))
\end{multline*}
be a sphere of type $(k,\ell)$.
We say that $\Sigma_1$ and $\Sigma_2$ are \emph{conformally equivalent} if
$m=k$, $n=\ell$ and there is
a complex analytic isomorphism $F:S_1 \to S_2$ such that $F(p_i)=q_i$ and
$\Phi_i = \Psi_i \circ F$ in some neighborhood of $p_i$ for $i=-m, \ldots, -1,1,\ldots,n$.

\subsection{Sewing spheres with tubes} \label{S:sewingspheres}
A fundamental property of the interactions of strings in conformal field theory is that certain
interactions should compose naturally.  In geometry we model composition by combining two
Riemann surfaces into a single Riemann surface.  This procedure is rigorously established
in \cite{H} for ``sewing", or attaching, two spheres of type $(1,m)$ and $(1,n)$ respectively.
However, Huang's description never uses the orientation on the punctures that remain
unsown, allowing the argument to be directly generalized to sewing one incoming
tube of a sphere of type $(m,n)$ to one outgoing tube of a sphere of type $(k, \ell)$ for
$m,\ell \in \N, \ n,k \in \Z_+$.
Following \cite{H}, we will describe the conditions under which such a
sewing may occur and describe the resulting sphere with tubes.

We will use $B^r \subset \C$ and $\overline{B}^r \subset \C$ to
denote the open and closed discs of radius $r$ centered at the
origin.  With $\Sigma_1$ and $\Sigma_2$ as above, choose integers
$1 \leq i \leq n$ and $1 \leq j \leq k$.  We say that the
\emph{$i$-th tube of $\Sigma_1$ can be sewn with the $-j$-th tube
of $\Sigma_2$} if there exists $r \in \R_+$ such that

\begin{align*}
\overline{B}^r &\subset \Phi_i(U_i), \\
\overline{B}^{1/r} &\subset \Psi_{-j}(V_{-j}),
\end{align*}

\noindent $p_i$ is the only puncture in
$\Phi_i^{-1}(\overline{B}^r)$, and $q_{-j}$ is the only puncture
in $\Psi_{-j}^{-1}(\overline{B}^{1/r})$. Additionally, we say the
$i$-th tube of $\Sigma_1$ can be sewn with the $-j$-th tube of
$\Sigma_2$ if there exists $r \in \R_+$ such that the domain of
$\Phi_i^{-1}$ may be analytically extended to $\overline{B}^r$
without $\Phi_i^{-1}(\overline{B}^r)$ containing punctures other
than $p_i$, and the domain of $\Psi_{-j}^{-1}$ may be analytically
extended to $\overline{B}^{1/r}$ without
$\Psi_{-j}^{-1}(\overline{B}^{1/r})$ containing punctures other
that $q_{-j}$. From such a $\Sigma_1$ and $\Sigma_2$ we obtain a
sphere with tubes of type $(m+k-1,n+\ell -1)$ as follows.  First,
choose $r_1, \ r_2$ real numbers such that $0 \ < \ r_2 \ < \ r \
< \ r_1$, ($\Phi_i$ can be extended so that) $\overline{B}^{r_1}
\subset \Phi_i(U_i)$, ($\Psi_{-j}$ can be extended so that)
$\overline{B}^{1/r_2} \subset \Psi_{-j}(V_{-j})$,
$\Phi_i^{-1}(\overline{B}^{r_1})$ contains no punctures of
$\Sigma_1$ besides $p_i$, and
$\Psi_{-j}^{-1}(\overline{B}^{1/r_2})$ contains no punctures of
$\Sigma_2$ besides $q_{-j}$.  Then define the sphere

\begin{equation*}
S_3=((S_1 \smallsetminus \Phi_i^{-1}(\overline{B}^{r_2}) \sqcup
(S_2 \smallsetminus \Psi_{-j}^{-1}(\overline{B}^{1/r_1})) \slash \sim
\end{equation*}

\noindent where (using the notation of \cite{V}) $\sqcup$ is the
disjoint union and $\sim$ is the equivalence relation, preserving
complex structure, given by $p \sim q$ if and only if $p=q$ or $p
\in \Phi_i^{-1}(\overline{B}^{r_1}) \smallsetminus
\Phi_i^{-1}(\overline{B}^{r_2})$, $q \in
\Psi_{-j}^{-1}(\overline{B}^{1/r_2}) \smallsetminus
\Psi_{-j}^{-1}(\overline{B}^{1/r_1})$ and
$\Psi_{-j}^{-1}(\frac{1}{\Phi_i(p)})=q$.  We define the sewing of
the $i$-th tube of $\Sigma_1$ with the $-j$-th tube of $\Sigma_2$
to be the sphere $S_3$ with the ordered punctures

\begin{equation*}
q_{-k}, \ldots, q_{-j-1}, p_{-m}, \ldots, p_{-1}, q_{-j+1}, \ldots, q_{-1},p_1, \ldots, p_{i-1},
q_1, \ldots, q_{\ell},p_{i+1}, \ldots, p_n
\end{equation*}

\noindent
and the local coordinate maps restricted appropriately; that is, for each puncture $p$ of $\Sigma_1$ remove
$\Phi_i^{-1}(\overline{B}^r_1)$ from the local coordinate neighborhood and restrict the local
coordinate map at $p$ to this new local coordinate neighborhood, and do similarly
for each puncture $q$ in $\Sigma_2$ and its local coordinate chart.  This sewing is independent of
the choice of $r, r_1, r_2$ (\cite{H}).

\subsection{The moduli space of spheres with tubes} \label{S:modulispace}
At first glance, spheres with tubes might appear to be exactly the
right picture for modelling closed strings in space-time.
However, a standard complex manifold structure turns out to be a
bit too rigid for the suggested interactions in space-time.  By
definition, conformal field theories assume invariance of
interactions under conformal transformations.  Thus two spheres
which possess an invertible conformal map between them will result
in equivalent correlation functions for the resulting particle
interaction.

First, we follow Huang's work (\cite{H2}, \cite{H}) describing canonical
representatives for each conformal equivalence class of spheres of type $(1,n)$, for $n \in \Z_+$.

For any $r \in \mathbb{R}_+$ and $z \in \C^{\times}$, let
\begin{align*}
&B_0^r=\{w \in \widehat{\C} | \ |w| < r\}, \\
&B_z^r=\{w \in \widehat{\C} | \ |z-w| < r\}, \\
&B_{\infty}^r=\{w \in \widehat{\C} | \ |1/w| < r\}.
\end{align*}

Recall that the reason we refer to
any compact genus-zero Riemann surface as a sphere is that each one is complex analytically isomorphic
to the standard sphere $\widehat{\C}= \C \cup \{ \infty \}$.  Furthermore, we have the following
proposition.

\begin{proposition} \label{P:canon1a}
Any sphere with tubes of type $(1,n)$, for $n \in \Z_+$, is conformally equivalent to a sphere of the form

\begin{equation} \label{E:sphere-midway}
(\widehat{\C} ; z_{-1} ,z_1, \ldots, z_{n}; (B_{z_{-1}}^{r_{-1}} ,\Phi_{-1}),
(B_{z_{1}}^{r_{1}},\Phi_{1}),
(B_{z_{2}}^{r_{2}},\Phi_{2}), \ldots, (B_{z_{n}}^{r_{n}},\Phi_{n}))
\end{equation}

\noindent
where $z_{-1}=\infty$, $z_{n} = 0$, ${z_i \in \C^{\times}}$ for $i=1,\ldots,n-1$,
satisfying $z_i \neq z_j$
for $i \neq j$, $r_i \in \mathbb{R}_+$ for $i=-1,1,2,\ldots,n$ and
$\Phi_{-1}, \Phi_{1}, \Phi_{2}, \ldots, \Phi_{n}$ are analytic on $B_{z_{-1}}^{r_{-1}},
B_{z_{1}}^{r_{1}}, B_{z_{2}}^{r_{2}}, \ldots, B_{z_{n}}^{r_{n}}$, respectively,
such that

\begin{align}
\Phi_{i}(z_i)=0 & \ \ \  i=-1, 1, \ldots , n \\
\lim_{w \to \infty} w\Phi_{-1}(w) =1, \label{E:a_triv1}\\
\lim_{w \to z_i} \frac{\Phi_i(w)}{w-z_i} \neq 0, & \ \ \ i=1, \ldots ,n.
\end{align}
\end{proposition}

This proposition is exactly what is proven in Proposition 1.3.1 of \cite{H}.

\begin{remark}
In Huang's work, \cite{H}, the outgoing puncture is labelled the
0-th puncture while the incoming punctures are labelled the first
through $n$-th punctures.  In order to generalize to multiple
outgoing punctures, we refer to outgoing punctures with negative
indices.  For example, in Proposition \ref{P:canon1a} the lone
outgoing puncture is referred to as $z_{-1}$ and its local
coordinate chart is referred to as $(B_{z_{-1}}^{r_{-1}}
,\Phi_{-1})$ whereas these would be referred to as $z_0$ and
$(B_{\infty}^{r_{0}} ,\Phi_{0})$, respectively, in \cite{H}.  Our
new notation is required to consider multiple outgoing punctures
but also highlights the natural symmetry between incoming and
outgoing punctures.
\end{remark}

The fact that every conformal equivalence class
of spheres of type $(1,n)$, for $n \in \Z_+$, contains a sphere of the form (\ref{E:sphere-midway})
allows us to focus solely on the Riemann sphere,
$\widehat{\C}$, with tubes and still include every equivalence class.  The following
proposition goes a step farther and describes the precise amount of information
needed to specify a particular equivalence class.

\begin{proposition}\label{P:canon2a}
Let
\begin{multline*}
\Sigma_1=( \widehat{\C} ; \infty,z_{1}, \ldots, z_{n-1}, 0;
(B_{\infty}^{r_{-1}} ,\Phi_{-1}),
(B_{z_{1}}^{r_{1}},\Phi_{1}), \ldots, \\
(B_{z_{n-1}}^{r_{n-1}},\Phi_{n-1}),
(B_{0}^{r_{n}},\Phi_{n}))
\end{multline*}
and
\begin{multline*}
\Sigma_2=(\widehat{\C} ; \infty, \zeta_{1}, \ldots, \zeta _{n-1}, 0;
(B_{\infty}^{s_{-1}} ,\Psi_{-1}), (B_{\zeta _{1}}^{s_{1}},\Psi_{1}), \ldots \\
(B_{\zeta_{n-1}}^{s_{n-1}},\Psi_{n-1}), (B_{0}^{s_{-n}},\Psi_{-n}))
\end{multline*}

\noindent
be two spheres of type $(1,n)$, for $n \in \Z_+$.  Let $f_{-1}, f_{1}, \ldots,
f_{n}$ and $g_{-1}, g_{1}, \ldots, g_{n}$ be the series obtained by
expanding the analytic functions $\Phi_{-1}, \Phi_{1}, \ldots, \Phi_{n}$ and
$\Psi_{-1}, \Psi_{1}, \ldots, \Psi_{n}$ around $w= \infty, z_{1}, \ldots,
z_{n-1}, 0$ and $w=\infty, \zeta_{1}, \ldots, \zeta_{n-1}, 0$,
respectively.  The worldsheets $\Sigma_1$ and $\Sigma_2$ are conformally equivalent
if and only if $z_{i} = \zeta_{i}$ for
$i=1, \ldots, n-1$, and $f_{i}=g_{i}$ (as power series), for $i=-1,1, \ldots, n$.
\end{proposition}

This proposition along with its proof may by found in \cite{H} (Proposition 1.3.3).
To extend the above results to spheres of type $(1,0)$, we need the following two
propositions (Propositions 1.3.4 and 1.3.6 in \cite{H}, respectively).

\begin{proposition}\label{P:canon3a}
Any sphere with tubes of type $(0,1)$ is conformally equivalent to a sphere of the form

\begin{equation*}
(\widehat{\C} ; \infty; (B_{\infty}^{r_{-1}} ,\Phi_{-1}))
\end{equation*}

\noindent
where $r_1 \in \mathbb{R}_+$ and $\Phi_1$ is an analytic function on $B_{\infty}^{r_{-1}}$
that can be expanded as

\begin{equation} \label{E:a_triv2}
\Phi_1(w)=\frac{1}{w} + \sum_{j=2}^{\infty} a_j \left( \frac{1}{w} \right)^{j+1}
\end{equation}

\noindent
with each $a_j \in \C$.
\end{proposition}

\begin{proposition}\label{P:canon4a}
Two spheres of the form $Q_1=(\widehat{\C} ; \infty; (B_{\infty}^{r_{-1}} ,\Phi_{-1}))$
and $Q_2=(\widehat{\C} ; \infty; (B_{\infty}^{s_{-1}} ,\Psi_{-1}))$ are conformally
equivalent if and only if $f_{-1}=g_{-1}$ (as power series), where $f_{-1}$ and $g_{-1}$ are the
power series expansions of $\Phi_{-1}$ and $\Psi_{-1}$, respectively, around $w=\infty$.
\end{proposition}

These propositions allow us to canonically describe conformal equivalence classes of spheres of
type $(1,n)$ for $n \in \N$ as follows:

\begin{equation} \label{E:123b}
(z_{1}, \ldots, z_{n-1}; f_{-1},f_{1}, \ldots, f_{n}) \ \ \ \text{ if } n > 0
\end{equation}
\begin{equation}
(f_{-1}) \ \ \ \text{ if } n=0
\end{equation}

\noindent
where $z_{1}, \ldots, z_{n-1}$ are distinct nonzero complex numbers and
$f_{-1},f_{1},  \ldots, f_{n}$ are power series that are convergent in some positive
neighborhood. These tuples represent the conformal equivalence classes of the
worldsheets

\begin{multline*}
(\widehat{\C} ; \infty , z_{1}, \ldots, z_{n-1}, 0;
(B_{\infty}^{r_{-1}} ,f_{-1}),(B_{z_{1}}^{r_{1}}, f_{1}), \\
\ldots, (B_{z_{n-1}}^{r_{n-1}},f_{n-1}),
(B_0^{r_{n}},f_{n})) \ \ \ \text{ if } n > 0
\end{multline*}

\begin{equation*}
(\widehat{\C} ; \infty; (B_{\infty}^{r_{-1}} ,f_{-1})) \ \ \ \text{ if } n=0,
\end{equation*}

\noindent
respectively, where $r_{-1},r_{1},\ldots,r_{n}$ are appropriately chosen radii of convergence
so that the corresponding local coordinate maps are convergent within
$B_{\infty}^{r_{-1}}$, $B_{z_{1}}^{r_{1}}, \ldots$, $B_{z_{n-1}}^{r_{n-1}},$
$B_0^{r_n}$, respectively, and the local coordinate neighborhoods do not overlap.
Notice that the lack of specificity on the choice of the $r_i$ means that we are not choosing
a specific canonical representative of each equivalence class in the most specific sense.
This subtlety is usually suppressed and the above tuples are referred to as canonical representatives
because the germ of each analytic function is the only data affecting which equivalence class a
particular sphere with punctures belongs to.

We now turn our attention to spheres with outgoing tubes.  The
results for spheres of type $(n,1)$, for $n \in \Z_+$ follow
directly from the spheres with incoming tubes case since the
proofs of those propositions never refer to the orientations of
punctures. In the case of spheres of type $(0,1)$, however,
additional proof is necessary in order to normalize the lone
incoming puncture at 0 instead of at $\infty$ as is the case for
spheres of type $(1,0)$.  (We normalize the single incoming
puncture in the representative of conformal equivalence classes of
spheres of type $(0,1)$ to be at 0 in order to maintain
compatibility with the preexisting structure notation for GVOAs
defined by Huang in \cite{H}.)

\begin{proposition} \label{P:canon1}
Any sphere with tubes of type $(n,1)$, for $n \in \Z_+$ is conformally equivalent to a
sphere of the form

\begin{equation*}
(\widehat{\C} ; z_{-n}, \ldots, z_{-1},z_{1}; ,
(B_{z_{-n}}^{r_{-n}},\Phi_{-n}),
 \ldots, (B_{z_{-1}}^{r_{-1}},\Phi_{-1}), (B_{z_1}^{r_1} ,\Phi_1))
\end{equation*}

\noindent
where $z_{-n} = \infty$, $z_1=0$, ${z_i \in \C^{\times}}$ for $i=-n+1, \ldots,-1$,
satisfying $z_i \neq z_j$
for $i \neq j$, $r_i \in \R_+$ for $i=-n,\ldots,-1,1$ and
$\Phi_{-n}, \ldots, \Phi_{-1}, \Phi_{1}$ are analytic on $B_{z_{-n}}^{r_{-n}},
\ldots, B_{z_{-1}}^{r_{-1}}, B_{z_{1}}^{r_{1}}$, respectively,
such that

\begin{align}
\Phi_{i}(z_i)=0 & \ \ \  i=-n,\ldots ,-1, 1 \\
\lim_{w \to \infty} w\Phi_{-n}(w) =1, \label{E:a_triv3} \\
\lim_{w \to z_i} \frac{\Phi_i(w)}{w-z_i} \neq 0, & \ \ \ i=-n+1,\ldots ,-1,1.
\end{align}
\end{proposition}

\begin{proof}  The proof of Proposition 1.3.1 in \cite{H} suffices since it never uses
the orientation of the punctures.
\end{proof}

As before, this proposition narrows the choices for a canonical representatives
in each equivalence class of spheres of type $(n,1)$ and the following proposition
illuminates the exact minimum information required to reference a particular
equivalence class.

\begin{proposition}\label{P:canon2}
Let

\begin{multline*}
\Sigma_1=( \widehat{\C} ; \infty , z_{-n+1}, \ldots, z_{-1}, 0;
(B_{\infty}^{r_{-n}} ,\Phi_{-n}),
(B_{z_{-n+1}}^{r_{-n+1}},\Phi_{-n+1}), \ldots, \\
(B_{z_{-1}}^{r_{-1}},\Phi_{-1}),
(B_0^{r_{1}},\Phi_{1}))
\end{multline*}
and
\begin{multline*}
\Sigma_2=(\widehat{\C} ; \infty, \zeta_{-n+1}, \ldots, \zeta _{-1}, 0;
(B_{\infty}^{s_{-n}} ,\Psi_{-n}), (B_{\zeta _{-n+1}}^{s_{-n+1}},\Psi_{-n+1}), \ldots, \\
(B_{\zeta_{-1}}^{s_{-1}},\Psi_{-1}),
(B_0^{s_{1}},\Psi_{1}))
\end{multline*}

\noindent
be two spheres of type $(n,1)$, for $n \in \Z_+$.  Let $f_{-n}, \ldots,
f_{-1}, f_{1}$ and $g_{-n}, \ldots, g_{-1}, g_{1}$ be the series obtained by
expanding the analytic functions $\Phi_{-n}, \ldots, \Phi_{-1}, \Phi_{1}$ and
$\Psi_{-n}, \ldots, \Psi_{-1}, \Psi_{1}$ around $w= \infty, z_{-n+1}, \ldots, z_{-1}, 0$
and $w= \infty, \zeta_{-n+1}, \ldots, \zeta_{-1}, 0$,
respectively.  The worldsheets $\Sigma_1$ and $\Sigma_2$ are conformally equivalent
if and only if $z_{i} = \zeta_{i}$, for
$i=-n+1, \ldots, -1$ and $f_{i}=g_{i}$ (as power series), for $i=-n, \ldots, -1, 1$.
\end{proposition}

\begin{proof}  Follows from Proposition 1.3.3 in \cite{H}.
\end{proof}

These two propositions allow us to choose a canonical representative for each conformal
equivalence class of worldsheets of type $(n,1)$ for $n \in \Z_+$.  Here we switch
notation slightly for power series centered at points other than 0 and $\infty$.
Power series expanded thus far about 0, or
about any nonzero complex number, have been power series in $w$, while the power series
at $\infty$ had to be expanded in terms of $\frac{1}{w}$.  However, for the nonzero complex
punctures, expanding in terms of $\frac{1}{w}$ is also valid.  In describing representatives
for equivalence classes of spheres with outgoing tubes, we will use the latter convention for
all outgoing (i.e., nonzero) punctures.  Thus a canonical representative will be denoted

\begin{equation} \label{E:123a}
(\frac{1}{z_{-n+1}}, \ldots, \frac{1}{z_{-1}}; f_{-n}, \ldots, f_{-1}, f_{1}),
\end{equation}

\noindent
where $z_{-n+1}, \ldots, z_{-1}$ are distinct nonzero complex numbers,
$f_{-n} \ldots, f_{-1}$ are power series in $\frac{1}{w}$ centered at $z_{-n+1}, \ldots, z_{-1}$,
respectively, that are convergent in some positive neighborhood, and $f_1$ is a power series in
$w$ centered at 0 that is convergent in some positive
neighborhood.  This canonical representative will represent the equivalence class of worldsheets
containing the worldsheet

\begin{multline*}
(\widehat{\C} ; \infty, \frac{1}{z_{-n+1}}, \ldots, \frac{1}{z_{-1}},
0; (B_{\infty}^{r_{-n}} ,f_{-n}),(B_{\frac{1}{z_{-n+1}}}^{r_{-n+1}},
f_{-n+1}), \\
\ldots, (B_{\frac{1}{z_{-1}}}^{r_{-1}},f_{-1}),
(B_0^{r_{1}},f_{1})),
\end{multline*}

\noindent
where $r_{-n},\ldots, r_{-1}, r_{1}$ are appropriately chosen radii of convergence
so that each local coordinate map is convergent within its corresponding
local coordinate neighborhood.

It turns out to be useful to refer to nonzero punctures as
$\frac{1}{z}$, first because we will be applying the global
transformation $w \mapsto \frac{1}{w}$ to the Riemann spheres used
in \cite{H} in some sense (cf. Remark \ref{R:operad_iso}), second
because it makes composing multiple shifts to $\infty$ clearer,
and third because it simplifies the isomorphism between geometric
vertex operator coalgebras and vertex operator coalgebras at the
conclusion of this paper.

For spheres of type $(0,1)$ we use a similar approach to that used
for spheres of type $(1,0)$ (studying linear fractional
transformations of the sphere) as the following two propositions
reveal.  (For explicit proofs, see \cite{K}.)

\begin{proposition}\label{P:canon3}
Any sphere with tubes of type $(0,1)$ is conformally equivalent to a sphere of the form

\begin{equation*}
(\widehat{\C} ; 0; (B_0^{r_1} ,\Phi_1))
\end{equation*}

\noindent
where $r_1 \in \mathbb{R}_+$ and $\Phi_1$ is an analytic function on $B_0^{r_1}$
that can be expanded as

\begin{equation} \label{E:a_triv4}
\Phi_1(w)=w + \sum_{j=2}^{\infty} a_j w^{j+1},
\end{equation}

\noindent
for $a_j \in \C$.
\end{proposition}

Proceeding in the same way as we did in Proposition \ref{P:canon4a} for worldsheets of type
$(1,0)$, we now argue:

\begin{proposition}\label{P:canon4}
Two spheres of the form $Q_1=(\widehat{\C} ; 0; (B_0^{r_1} ,\Phi_1))$
and $Q_2=(\widehat{\C} ; 0; (B_0^{s_1} ,\Psi_1))$ are conformally equivalent
if and only if $f_1=g_1$ (as power series), where $f_{1}$ and $g_{1}$ are the
power series expansions of $\Phi_1$ and $\Psi_1$, respectively, around $w=0$.
\end{proposition}

A representative for an equivalence class of spheres of type $(0,1)$ will be described
as

\begin{equation*}
(f_{1})
\end{equation*}

\noindent
which will represent the equivalence class of spheres of type $(1,0)$ containing

\begin{equation*}
(\widehat{\C} ; 0; (B_0^{r_{1}} ,f_{1}))
\end{equation*}

\noindent
where $r_{1}$ is an appropriately chosen radius of convergence.

\begin{remark}
Although we will be primarily concerned with the moduli space of spheres with incoming tubes
and the moduli space of spheres with outgoing tubes separately, the above propositions allow
us to describe the general moduli space of spheres of type $(m,n)$, for any $m,n \in \N$.
The type $(0,0)$ case is a one element set (that is not sewable); the cases of spheres of types
$(1,0)$ and $(0,1)$ have been described above; and given $m,n \in \Z_+$, the moduli spaces of
spheres of type $(m,n)$ may be describes as

\begin{equation} \label{E:123}
(\frac{1}{z_{-m+1}}, \ldots, \frac{1}{z_{-1}}, z_1, \ldots, z_{n-1}; f_{-m}, \ldots, f_{-1}, f_{1},
\ldots, f_{n})
\end{equation}

\noindent
by making the trivial generalization to Proposition \ref{P:canon1a} or \ref{P:canon1}
that the orientations of punctures other than the first and last are irrelevant.  Here Equation
(\ref{E:123}) is the generalization of Equations (\ref{E:123b}) and (\ref{E:123a}), and should
be understood similarly.
\end{remark}

\subsection{Sewing moduli spaces of spheres with tubes} \label{S:sewingmoduli}

Now that we have described the moduli space of spheres with tubes
(the moduli space of spheres with incoming tubes and the moduli
space of spheres with outgoing tubes being subsets), we need to
lift the sewing operation from the sewing of particular spheres to
the sewing of elements of the moduli space. Given $\Sigma_1$ and
$\Sigma_2$, two worldsheets such that the $j$-th outgoing puncture
of $\Sigma_2$ may be sewn into the $i$-th incoming puncture of
$\Sigma_1$ as in Section \ref{S:sewingspheres}, let $\Sigma_3$ be
the worldsheet resulting from the sewing.  We need to be certain
that given $\Sigma'_1$ and $\Sigma'_2$ in the same conformal
equivalence classes as $\Sigma_1$ and $\Sigma_2$, respectively,
that the $j$-th outgoing puncture of $\Sigma'_2$ may be sewn into
the $i$-th incoming puncture of $\Sigma'_1$ and that the resulting
sphere, $\Sigma'_3$, is in the same equivalence class as
$\Sigma_3$.  But given complex analytic isomorphisms $G_1:
\Sigma_1 \to \Sigma'_1$ and $G_2: \Sigma_2 \to \Sigma'_2$
witnessing conformal equivalence, a radius $r \in \R_+$ that make
$\Sigma_1$ and $\Sigma_2$ sewable is exactly one that makes
$\Sigma'_1$ and $\Sigma'_2$ sewable.  Further, $G_1 \sqcup G_2$
will provide precisely the complex analytic isomorphism needed
between $\Sigma_3$ and $\Sigma'_3$.

When it is the case that two elements of the moduli space of spheres, $Q_1$ and $Q_2$, are sewable
with the $j$-th outgoing tube of $\Sigma_2$ attaching to the $i$-th incoming tube of $\Sigma_1$, we
denote the resulting sewing by $Q_1 \ _i \infty_{-j} \ Q_2$ (cf. \cite{V}, \cite{H}).
Now that we have established that this sewing operation, under appropriate conditions, is
well-defined on elements of the moduli space of
spheres with punctures, let

\begin{align*}
Q_1 &= (z_{-m+1}^{-1}, \ldots, z_{-1}^{-1}, z_1, \ldots, z_{n-1}; f_{-m}, \ldots, f_{-1},
f_{1}, \ldots, f_{n}) \\
Q_2 &= (\zeta_{-k+1}^{-1}, \ldots, \zeta_{-1}^{-1}, \zeta_1, \ldots, \zeta_{\ell-1}; g_{-k}, \ldots, g_{-1},
g_{1}, \ldots, g_{\ell})
\end{align*}

\noindent be two elements of the moduli space with spheres such
that the sewing $Q_1 \ _i \infty_{-j} Q_2$ exists.  Equivalently
we may require that given the canonical representatives of $Q_1$
and $Q_2$, for $\ell,m \in \N$, $k,n \in \Z_+$ and $1 \leq i \leq
n$, $1 \leq j \leq k$, there exists $r >0$ with
$f_{i}^{-1}(\overline{B}^{r})$ and
$g_{-j}^{-1}(\overline{B}^{1/r})$ well defined and containing only
the punctures $z_i$ and $\zeta_{-j}^{-1}$, respectively (note that
we consider $z_n=0$ and $\zeta_{-k}^{-1}=\infty$).  In this case,
there are $r_1$ and $r_2$ satisfying $0<r_2<r<r_1$ such that
$f_{i}^{-1}(\overline{B}^{r_1})$ and
$g_{-j}^{-1}(\overline{B}^{1/r_2})$ are well defined and still
contain only the punctures $z_i$ and $\zeta_{-j}^{-1}$,
respectively.  Choose

\begin{equation*}
Q_3= (\tilde{z}_{-m-k+2}^{-1}, \ldots, \tilde{z}_{-1}^{-1}, \tilde{z}_1, \ldots,
\tilde{z}_{n+\ell-2}; \tilde{f}_{-m-k+1}, \ldots, \tilde{f}_{-1},
\tilde{f}_{1}, \ldots, \tilde{f}_{n+\ell-1})
\end{equation*}

\noindent so that $Q_1 \ _i \infty_{-j} Q_2=Q_3$.  Let $F$ be the
unique conformal equivalence from the sewing of the canonical
representatives of $Q_1$ and $Q_2$ to the canonical representative
of $Q_3$. Since the sewing is in two components, there are also
unique maps $F^{(1)}:\widehat{\C} \setminus
f_i^{-1}(\overline{B^{r_2}}) \to \widehat{\C}$ and
$F^{(2)}:\widehat{\C} \setminus
g_{-j}^{-1}(\overline{B^{1/r_1}}\{\zeta_{-i}^{-1} \}) \to
\widehat{\C}$ satisfying

\begin{equation} \label{E:sewing_eqn1}
\left. F^{(1)}(w) \right|_{f_i^{-1}(B^{r_2} \setminus \overline{B}^{r_1})}=
\left. F^{(2)}\left(g_{-j}^{-1}\left(\frac{1}{f_i(w)}\right)\right)
\right|_{f_i^{-1}(B^{r_2} \setminus \overline{B}^{r_1})}
\end{equation}
\begin{equation} \label{E:normalize1}
F^{(1)}(\infty)=\infty
\end{equation}
\begin{equation} \label{E:normalize2}
F^{(2)}(0)=0
\end{equation}
\begin{equation} \label{E:normalize3}
\lim_{w \to \infty} \frac{F^{(1)}(w)}{w}=1.
\end{equation}

\noindent We call (\ref{E:sewing_eqn1}) the \emph{sewing equation}
and (\ref{E:normalize1}), (\ref{E:normalize2}),
(\ref{E:normalize3}) the \emph{normalization conditions}.  These
are the same equations as developed in Section 1.4 of \cite{H} but
in a slightly more general context. A good deal of work goes into
describing the canonical representative of $Q_{1} \
_{i}\infty_{-j} \ Q_2$ constructively, but this work is quite
instructive and indeed necessary in several cases for the main
isomorphism theorem between GVOCs and VOCs. We will include a pair
of example of explicit sewing in Section \ref{S:example}.

In nearly any structure with an operation, it is valuable to investigate
what set together with the available operations generates all possible elements.  That
question is answered for the moduli space of spheres with incoming tubes in Proposition 1.3.9
of \cite{H} as follows.

\begin{proposition} \label{P:gen_mod_space1}
Any element of the moduli space of spheres of type $(1,n)$, for $n \in \N$, is generated under sewing
by the
element $(w^{-1})$ of type $(1,0)$, elements of type $(1,1)$, and elements of type $(1,2)$ which are of the
form $(z; w^{-1}, w-z, w)$, for $z \in \C^{\times}$.
\end{proposition}

In the proposition $w^{-1}$, $w-z$, and $w$ should be understood as power series in $w^{-1}$, $w$,
and $w$ centered at $\infty$, $z$ and 0, respectively.
We are able to demonstrate a similar generating set for the moduli space of spheres with outgoing
tubes (using similar notation for power series).

\begin{proposition} \label{P:gen_mod_space}
Any element $Q$ of the moduli space of spheres of type $(n,1)$, for $n \in \N$, is generated under
sewing by the element $(w)$ of type $(0,1)$, elements of type $(1,1)$, and elements
of type $(2,1)$ which are of the form
$(z^{-1}; w^{-1}, w^{-1}-z, w)$, for $z \in \C^{\times}$.
\end{proposition}

\begin{proof}
The proof is analogous to the proof of Proposition 1.3.9 in
\cite{H}, i.e., by using a Jordan curve to split a given sphere
(cf. \cite{K}).
\end{proof}

\subsection{Expansions in terms of infinitesimal local coordinates}
Any local coordinate map at 0 may be expressed uniquely as

\begin{equation} \label{E:exp1}
\exp\left(  \sum_{j \in \Z_+} A_j w^{j+1} \frac{d}{dw}\right)
a_0^{w\frac{d}{dw}} w
\end{equation}

\noindent
where  the $A_j \in \C$ and $a_0 \in \C^{\times}$.
(See Proposition 2.1.1 in \cite{H} and the following discussion).

Any local coordinate map vanishing at $\infty$, on the other hand, may be written as

\begin{equation} \label{E:exp2}
\exp\left(-  \sum_{j \in \Z_+} A_j w^{-j+1} \frac{d}{dw}\right)
 a_0^{-w\frac{d}{dw}} \frac{1}{w}
\end{equation}

\noindent
where  the $A_j \in \C$ and $a_0 \in \C^{\times}$.
(See Proposition 2.1.16 in \cite{H}).  In Huang's treatment of these moduli
spaces as well as in this work, the puncture at $\infty$ is usually normalized
so that $a_0 = 1$.  This is a direct consequence of the conditions (\ref{E:a_triv1}) and
(\ref{E:a_triv2}) in the moduli space of spheres with incoming tubes and the condition
(\ref{E:a_triv3}) in the moduli space with outgoing tubes.

Any local coordinate map vanishing at a nonzero complex number, $z$, may be written as a local
coordinate map at 0 composed with a shift of $z$ to 0:

\begin{equation*}
\left. \exp\left(  \sum_{j \in \Z_+} A_j x^{j+1} \frac{d}{dx}\right)
 a_0^{x \frac{d}{dx}} x \right|_{x= w - z},
\end{equation*}

\noindent
or as a local
coordinate map at $\infty$ composed with a shift of $z$ to $\infty$:

\begin{equation*}
\left. \exp\left(-  \sum_{j \in \Z_+} A_j x^{-j+1} \frac{d}{dx}\right)
 a_0^{-x \frac{d}{dx}} \frac{1}{x} \right|_{x= (w^{-1} - z^{-1})^{-1}}.
\end{equation*}

We will use the expansion at 0 for incoming tubes (as in \cite{H})
and the expansion at $\infty$ for outgoing tubes.  Let the
sequence $A= \{A_j\}_{j\in\Z_+}$ denote the higher order
coefficients in a given exponential expansions in terms of the
infinitesimal local coordinates, $x^{j+1} \frac{d}{dx}$, and let
the notation $(a_0,(A_1,A_2,\ldots))$ or $(a_0,A)$ record all the
coefficients in a given expansion.  We will also use the notation
$\mathbf{0}$ for the sequence of all zeros and $ A (a)$ for $\{
a^j A_j \}_{j\in\Z_+}$.  When we are dealing with a sequence of
formal variables, we will denote the sequence
$(\alpha_0,(\mathcal{A}_1,\mathcal{A}_2,\ldots))$ or
$(\alpha_0,\mathcal{A})$ following \cite{H} and \cite{B}.

Using this notation, the canonical representative for an equivalence class of genus-zero
worldsheets of type $(1,n)$ can be expanded as

\begin{equation} \label{E:canon1a}
(z_{1}, \ldots, z_{n-1}; A^{(-1)}, (a_0^{(1)},A^{(1)}), \ldots,
(a_0^{(n)},A^{(n)})),
\end{equation}

\noindent
and the canonical representative for an equivalence class of type $(n,1)$ can be expanded as

\begin{equation} \label{E:canon1}
(z_{-n+1}^{-1}, \ldots, z_{-1}^{-1}; A^{(-n)},
(a_0^{(-n+1)},A^{(-n+1)}), \ldots, (a_0^{(-1)},A^{(-1)}),
(a_0^{(1)},A^{(1)})),
\end{equation}

\noindent with $n\in\Z_+$, the $z_{i}$'s distinct nonzero complex
numbers, $(a_0^{i},A^{(i)})$ recording the local coordinate map
for the puncture of corresponding index, and for the first
sequence in each list, $A^{(-1)}$ and $A^{(-n)}$ respectively,
recording the higher order coefficients for the exponential
expansion of the local coordinate map at $\infty$.  As mentioned
before, these local coordinate maps will always be normalized to
have $a_{0}^{(-1)}=1$ and $a_0^{(-n)}=1$, respectively.
Exponential notation also allows us to describe a canonical
representative of type $(1,0)$ as

\begin{equation} \label{E:canon2a}
(A^{(-1)})
\end{equation}

\noindent where $A^{(-1)}$ records the coefficients of the
exponential expansion at $\infty$, and to describe a canonical
representative of type $(0,1)$ as

\begin{equation} \label{E:canon2}
((1,A^{(1)})),
\end{equation}

\noindent where $(1,A^{(1)})$ records the local coordinate map for
the puncture at 0. We see from Proposition \ref{P:canon3a} that in
the type $(1,0)$ case not only must
 $a_{0}^{(-1)}=1$ but $A_1^{(-1)}=0$ for all representatives.  Also, Proposition \ref{P:canon3}
implies that in the type $(0,1)$ case $a_0^{(1)}=1$ and $A_1^{(1)}=0$.

\subsection{A formal description of the moduli space of spheres with incoming and outgoing tubes}

We begin by formalizing a notation for describing the
moduli space of spheres with one outgoing and $n$ incoming tubes. For $n \in \Z_+$, let

\begin{equation*}
M^{n-1} = \{(z_{1},z_{2}, \ldots z_{n-1}) | z_i \in \C^{\times},
z_i \neq z_j \text{ for } i \neq j \}
\end{equation*}

\noindent
Note that for $n=1$, $M^0$ has exactly one element.  We say that the series \\
$\exp\left(  \sum_{j \in \Z_+} A_j w^{j+1}\frac{d}{dw}\right)w$ is absolutely convergent
in a neighborhood of 0 precisely if its expansion as powers of $w$ is absolutely
convergent in a neighborhood of 0.  We then define

\begin{align*}
H= \text{\Huge{\{}} A=\{A_j\}_{j \in \Z_+} \in \C^{\infty} \mid
\exp\left(  \sum_{j \in \Z_+} A_j w^{j+1}
\frac{d}{dw}\right)w \text{ is absolutely} \\
\text{ convergent in some neighborhood of } 0 \text{\Huge{\}}}.
\end{align*}

\begin{proposition}
The moduli space of spheres with tubes of type $(1,n)$, for $n \in \Z_+$,
can be identified with the set

\begin{equation*}
K(n)=M^{n-1} \times H \times (\C^{\times} \times H)^n,
\end{equation*}

\noindent
and the moduli space of spheres with tubes of type $(1,0)$ can be identified with
the set

\begin{equation*}
K(0)= \{ A \in H | A_1=0\}.
\end{equation*}
\end{proposition}

\noindent
This is simply a restatement of Propositions 3.1.1 and 3.1.2 in \cite{H} but is the obvious
conclusion of Equations (\ref{E:canon1a}) and (\ref{E:canon2a}).  Let

\begin{equation*}
K= \coprod_{n \in \N} K(n)
\end{equation*}

\noindent
be the moduli space of incoming tubes with spheres.  Under the sewing operation (which corresponds
to the `$\circ_i$' notation of operads), $K$ is a partial operad (\cite{HL} Section 5).
The action of $\sigma \in S_n$ on $S$, a genus-zero
Riemann surface with one outgoing and $n$ incoming punctures, is given by reordering the punctures and
making the $i$-th incoming puncture the $\sigma(i)$-th incoming puncture for each $i=1, \ldots,n$.  This action is invariant under
conformal equivalence and, hence, induces an action of $\sigma$ on the moduli space of spheres of type
$(1,n)$, that is $K(n)$.

 Via Equations (\ref{E:canon1}) and (\ref{E:canon2}),
we also establish a uniform way of describing a canonical
representative of the moduli space of spheres with one incoming
and $n$ ordered outgoing punctures.

\begin{proposition}
Let $\widetilde{M}^{n-1} = \{(\frac{1}{z_{-n+1}}, \ldots
\frac{1}{z_{-1}}) | z_i \in \C^{\times}, z_i \neq z_j \text{ for }
i \neq j \}$ and again let $H= \{A=\{A_j\}_{j \in \Z_+} |
\exp\left(  \sum_{j \in \Z_+} A_j w^{j+1} \frac{d}{dw}\right)w$ is
absolutely convergent in some neighborhood of 0 \}. Then the
moduli space of spheres with tubes of type $(n,1)$, for $n \in
\Z_+$, can be identified with the set

\begin{equation*}
K^*(n)=\widetilde{M}^{n-1} \times H \times (\C^{\times} \times H)^n,
\end{equation*}

\noindent
and the moduli space of spheres with tubes of type $(0,1)$ can be identified with
the set

\begin{equation*}
K^*(0)= \{A \in H | A_1=0\}.
\end{equation*}
\end{proposition}

This is simply the obvious interpretation of Equations (\ref{E:canon1}) and (\ref{E:canon2}).
As in the case of $K$, we define

\begin{equation*}
K^*= \coprod_{n \in \N} K^*(n)
\end{equation*}

\noindent
to be the moduli space of spheres with outgoing tubes.  Again, the sewing operation gives
$K^*$ the structure of a partial operad.  The action of $\sigma \in S_n$ on $K^*(n)$
is the same as that of $K$, induced from
making the $i$-th outgoing puncture the $\sigma(i)$-th outgoing puncture for each $i=1, \ldots,n$
on any representative of the equivalence class acted upon.
This action is invariant under conformal equivalence and, hence, induces an action of $\sigma$ on
$K^*(n)$.
See \cite{K} for a more detailed discussion.

The next two propositions follow directly from observations on spheres with tubes
and the fact that sewing and
permutation are invariant under conformal equivalence.  (They may also be thought of as extensions of
the results of Huang on $K$ to $K^*$ by observing that orientation of punctures is not at issue.)

\begin{proposition} \label{P:K^*_assoc}
Let $\ell,m \in \N$ and $n \in \Z_+$ such that $m+n > 1$.  Choose
$Q_1 \in K^*(\ell), \ Q_2 \in K^*(m), \ Q_3 \in K^*(n)$ and integers $1 \leq i \leq m+n-1$
and $1 \leq j \leq n$.  Then $Q_1 \ _1\infty_{-i} \ (Q_2 \ _1\infty_{-j} Q_3)$
exists if and only if one of the following 3 holds:

\noindent
(1) $i < j$, $Q_2 \ _1\infty_{-j-\ell+1} \ (Q_1 \ _1\infty_{-i} Q_3)$ exists and

\begin{equation*}
Q_1 \ _1\infty_{-i} \ (Q_2 \ _1\infty_{-j} Q_3)
=Q_2 \ _1\infty_{-j-\ell+1} \ (Q_1 \ _1\infty_{-i} Q_3);
\end{equation*}

\noindent
(2) $j \leq i < j+m$, $(Q_1 \ _1\infty_{-i+j-1} \ Q_2) \ _1\infty_{-j} Q_3$
exists and

\begin{equation*}
Q_1 \ _1\infty_{-i} \ (Q_2 \ _1\infty_{-j} Q_3)
=(Q_1 \ _1\infty_{-i+j-1} \ Q_2) \ _1\infty_{-j} Q_3;
\end{equation*}

\noindent
(3) $i \geq j+m$, $Q_2 \ _1\infty_{-j} \ (Q_1 \ _1\infty_{-i+m-1} Q_3)$
exists and

\begin{equation*}
Q_1 \ _1\infty_{-i} \ (Q_2 \ _1\infty_{-j} Q_3)
=Q_2 \ _1\infty_{-j} \ (Q_1 \ _1\infty_{-i+m-1} Q_3).
\end{equation*}
\end{proposition}

\begin{proposition} \label{P:K^*perm}
Let $Q_1, Q_2 \in K^*(2)$, $\sigma = (12) \in S_2$, and $\tau
=(132) \in S_3$. Then $Q_{1} \ _1\infty_{-1} \ Q_2$ exists if and
only if $\tau(Q_{1} \ _1\infty_{-2} \ (\sigma \ Q_2))$ also
exists.  If this is the case,

\begin{equation*}
Q_{1} \ _1\infty_{-1} \ Q_2= \tau(Q_{1} \ _1\infty_{-2} \ (\sigma
\ Q_2)).
\end{equation*}
\end{proposition}

Proposition \ref{P:K^*perm} may be extended to all of the
symmetric groups and spheres with outgoing tubes since every
symmetric group is generated by transpositions and every sphere
with outgoing tubes is generated as in Proposition
\ref{P:gen_mod_space}.

\begin{remark} \label{R:operad_iso}
It is significant to note that $K$ and $K^*$ are isomorphic as
partial operads.
%The isomorphism is induced by reversing the
%orientation of each puncture on worldsheets.  This is conformally
%invariant and, in fact, given $S_1$ and $S_2$, two spheres with
%incoming tubes, the conditions for the sewing $S_1 \ _i\infty_{-1}
%\ S_2$ to exist are precisely the same conditions required for the
%sewing $I(S_2) \ _{1}\infty_{-i} \ I(S_1)$ to exist with
The isomorphism is given by the global transformation
$I:\widehat{\C} \to \widehat{\C}:$ $w \mapsto \frac{1}{w}$ on the
underlying Riemann spheres of the canonical representatives of $K$
and $K^*$; in addition the isomorphism reverses orientation of
punctures, changes each puncture location $z$ to $\frac{1}{z}$,
and composes local coordinate maps with $I$.  Finally, given $Q
\in K(n)$, with local coordinates at 0 given by
$(a_0^{(n)},A^{(n)})$, renormalization via the map
$(a_0^{(n)})^{-w \frac{d}{dw}}$ is required to express $I(Q)$ in
its canonical form, but the isomorphism is clearest when canonical
representatives are not chosen.  The transformation map is
functorial on equivalence classes, i.e. it commutes with actions
of the symmetric groups by construction and also commutes with the
sewing operation.  To see that $I$ indeed commutes with sewing,
note that given $S_1$ and $S_2$, two spheres with incoming tubes,
the conditions for the sewing $S_1 \ _i\infty_{-1} \ S_2$ to exist
are precisely the same conditions required for the sewing $I(S_2)
\ _{1}\infty_{-i} \ I(S_1)$ to exist and, in fact, when both exist

\begin{equation*}
I(S_1 \ _i\infty_{-1} \ S_2)=I(S_2) \ _{1}\infty_{-i} \ I(S_1)
\end{equation*}

It is also useful to recall that $K^*(1)$ and $K(1)$ are both
names for equivalence classes under conformal isomorphism of
spheres of type $(1,1)$, i.e. $K^*(1)$ and $K(1)$ are equal as
sets with identical action of the Virasoro algebra, and identical
compositions $K^*(1) \times K^*(1) \to K^*(1)$ and $K(1) \times
K(1) \to K(1)$. (But note that $I$ is not the identity map between
$K(1)$ and $K^*(1)$.)
%
%While this isomorphism points out that the conformal geometry in
%this paper is strikingly similar to that found in \cite{H}, the
%difference in operadic structures in the formal algebra is what
%makes this work truly distinct from \cite{H}.  In \cite{H}, the
%fundamental building block of a geometric vertex operator algebra
%is a map from $\nu: K(n) \to Hom(V^{\otimes n}, \overline{V})$.
%Here, the fundamental building block will be a map $\mu: K^*(n)
%\to Hom(V, \overline{V^{\otimes n}})$
\end{remark}

\subsection{Commutation in the Virasoro algebra}
The main obstacle to describing the local coordinate maps of
$Q_{1} \ _1\infty_{-i} \ Q_2$  in exponential notation is the fact
that the operators $\{ w^{j+1}\frac{d}{dw} \}_{j \in \Z}$ do not
commute.  This was overcome by a somewhat specialized approach in
\cite{H} but since that time a more powerful method has been
introduced in \cite{BHL} which applies to the
Baker-Cambell-Hausdorff formula in general. This result is
applicable because the operators $\{-w^{j+1}\frac{d}{dw}\}_{j \in
\Z} $ give a representation of the Virasoro algebra.

Let $\mathcal{V}= \oplus_{n \in \Z} \C L_n \oplus \C d$ be the
Virasoro algebra with the usual commutation relations

\begin{align*}
[L_m, L_n] &= (m-n)L_{m+n} + \frac{m^3-m}{12} \delta_{m,-n}d, \\
[\mathcal{V},d] &=0
\end{align*}

\noindent for $m,n \in \Z$.  For commuting formal variables
$\{\mathcal{A}_k \}_{k \in \Z_+}$, $\{\mathcal{B}_k\}_{k \in
\Z_+}$, $\alpha_0$ and $\beta_0$, by \cite{BHL} there exist unique
$\{\Psi_k\}_{k \in \Z}$ and $\Gamma$ such that

\begin{multline} \label{E:4162}
e^{-\sum_{j \in Z_+} \mathcal{A}_j L_j} \alpha_0^{-L_0}
\beta_0^{-L_0} e^{-\sum_{j \in \Z_+} \mathcal{B}_j L_{-j}}  \\
=e^{-\sum_{j \in Z_+} \Psi_{-j} L_{-j}}(\alpha_0 \beta_0)^{L_0}
e^{\Psi_0 L_0} e^{-\sum_{j \in Z_+} \Psi_j L_j} e^{\Gamma d}.
\end{multline}

\noindent Given $Q_1$ and $Q_2$ such that an incoming puncture of
$Q_1$ with local coordinates $(a_0, A)$ may be sewn with an
outgoing puncture of $Q_2$ having local coordinates $(b_0, B)$,
these local coordinates may be substituted into the right-hand
side of (\ref{E:4162}) and the corresponding $\{\Psi_k\}_{k \in
\Z}$ and $\Gamma$ converge.

\begin{remark}
In \cite{H}, the convergence of $\{\Psi_k\}_{k \in \Z}$ and $\Gamma$ are shown for
sewings of $K(1)$ and then extend naturally to all sewings in $K$.  But $K(1)$ and $K^*(1)$
are one and the same (i.e. conformal equivalence classes of spheres of type $(1,1)$), so
convergence of $\{\Psi_k\}_{k \in \Z}$ and $\Gamma$ for sewings of $K^*(1)$ is already shown
and extends naturally to all of $K^*$. For the proof of convergence in $K(1)$, see
Corollary 4.3.2, Lemma 5.2.1 and p. 123 in \cite{H}.
\end{remark}

Because $\Gamma$ depends on $\mathcal{A}$, $\mathcal{B}$,
$\alpha_0$, and $\beta_0$, we denote $\Gamma$ by
$\Gamma(\mathcal{A}, \mathcal{B}, \alpha_0 \beta_0)$ following the
notational convention of \cite{H}. The definition of GVOA, as well
as GVOC, may be interpreted as depending on the choice of a
section of the determinant line bundle over each equivalence class
of spheres with tubes, in which case $d$ may be represented by the
complex value that determines the section and $L_k$ is represented
by $-w^{k+1}\frac{d}{dw}$ for $k \in \Z$.

\subsection{Examples of explicit sewing and the sewing identities}
\label{S:example}

It is instructive to sew a few simple spheres together in order to
observe the explicit local coordinates that are generated.

\begin{example} \label{E:sewing4}
Let

\begin{equation*}
Q_1=( \mathbf{0}, (a_0^{(1)}, A^{(1)}) ) =( w^{-1}, f_1) \in
K^*(1),
\end{equation*}
\begin{equation*}
Q_2=(\zeta^{-1};\mathbf{0}, (1,\mathbf{0}), (1,\mathbf{0})  )
=(\zeta^{-1};w^{-1}, w^{-1}-\zeta, w) \in K^*(2),
\end{equation*}

\noindent such that $Q_1 \ _1\infty_{-2} \ Q_2$ exists. The unique
$F^{(1)}$ and $F^{(2)}$ satisfying (\ref{E:sewing_eqn1}) are

\begin{align*}
F^{(1)}(w)=&w \\
F^{(2)}(w)=& f_1^{-1}(w)= (a_0^{(1)})^{-w \frac{d}{dw}}
e^{-\sum_{j \in \infty} A_j^{(1)} w^{j+1} \frac{d}{dw}}w.
\end{align*}

\noindent Thus

\begin{equation*}
Q_1 \ _1\infty_{-2} \ Q_2 = (f_1^{-1}(\zeta^{-1});\mathbf{0},
((a_0^{(1)})^{-1} e^{\Theta_0^{(2)}},
\Theta^{(2)}((a_0^{(1)})^{-1}) ), (a_0^{(1)}, A^{(1)}) )
\end{equation*}

\noindent where $\Theta_0^{(2)}$ and $\Theta^{(2)}=\{
\Theta_j^{(2)} \}_{j \in \Z_+}$ are determined by

\begin{multline*}
\left. e^{-\sum_{j \in \Z_+} (a_0^{(1)})^{-j} \Theta_j^{(2)}
x^{-j+1} \frac{d}{dx}} (a_0^{(1)})^{x \frac{d}{dx}}
e^{-\Theta_0^{(2)} x \frac{d}{dx}}\frac{1}{x}
\right|_{\frac{1}{x}=\frac{1}{w}-\frac{1}{f_{1}^{-1}(\zeta)}}= \\
e^{\sum_{j \in \Z_+} A_j^{(1)} w^{j+1} \frac{d}{dw}}\frac{1}{w}
-\zeta.
\end{multline*}

\noindent Equivalently, let $\hat{f}_2(x)=e^{\sum_{j \in \Z}
(a_0^{(1)})^{-j} A_j^{(1)} x^{-j+1} \frac{d}{dx}}x$ so that

\begin{equation} \label{E:ex4}
Q_1 \ _1\infty_{-2} \ Q_2 = ((a_0^{(1)}\hat{f}_2(\zeta))^{-1};
\mathbf{0}, ((a_0^{(1)})^{-1} e^{\Theta_0^{(2)}},
\Theta^{(2)}((a_0^{(1)})^{-1}) ), (a_0^{(1)},A^{(1)}) )
\end{equation}

\noindent where $\Theta_0^{(2)}$ and $\Theta^{(2)}=\{
\Theta_j^{(2)} \}_{j \in \Z_+}$ are determined by

\begin{equation*}
\left. e^{-\sum_{j \in \Z_+} \Theta_j^{(2)} x^{-j+1} \frac{d}{dx}}
e^{-\Theta_0^{(2)} x \frac{d}{dx}}\frac{1}{x}
\right|_{\frac{1}{x}=\frac{1}{w}-\hat{f}_2(\zeta^{-1})}=
(\hat{f}_2(w))^{-1} -\zeta.
\end{equation*}

%Note that this equation is simply a rewriting of Equation (2.2.30)
%in \cite{H} from the GVOA case.
Further, when considering formal
variables $\alpha_0^{(1)}$, $\mathcal{A}^{(1)}$, and $\zeta$
instead of $a_0^{(1)}$, $A^{(1)}$, and $\zeta$, respectively, the
$\Theta_j^{(2)}$'s satisfy the sewing identity called the second
sewing identity:

\begin{proposition} \label{P:sew_id3}
In the algebra $(\text{End } \C)[w,w^{-1}][\zeta] [\alpha_0^{(1)},
\alpha_0^{(1)}] [[\mathcal{A}^{(1)}]]$, we have

\begin{multline} \label{E:secondsewingidentity}
e^{-\sum_{k=-1}^{\infty} \left( \sum_{j \in \Z_+}
(\alpha_0^{(1)})^{-j} \mathcal{A}^{(1)}_j
{ -j+1 \choose k+1}\zeta^{-j-k} \right)w^{-k+1} \frac{\partial}{\partial w}} \\
=e^{ \Theta^{(2)}_0 w\frac{\partial}{\partial w}} e^{\sum_{j \in
\Z_+} \Theta_j^{(2)} w^{-j+1} \frac{\partial}{\partial w}}
e^{(\zeta-\hat{f}_2(\zeta))w^2 \frac{\partial}{\partial w}} .
\end{multline}
\end{proposition}

\begin{proof}
This identity is observed from Proposition 2.2.9 of \cite{H} as
follows. Let $x=w^{-1}, \ y=\zeta$, and
$\mathcal{B}=\mathcal{A}^{(1)}((\alpha_0^{(1)})^{-1})$ in
Proposition 2.2.9, then take the inverse of both sides of the
equation.
\end{proof}
\end{example}

Using the second sewing equation, (\ref{E:secondsewingidentity}),
we can generalize our example.

\begin{example} \label{E:sewing5}
Let

\begin{equation*}
Q_1=( A^{(-1)}, (a_0^{(1)},A^{(1)}) ) \in K^*(1),
\end{equation*}
\begin{equation*}
Q_2=(\zeta^{-1}; \mathbf{0}, (b_0^{(-1)},B^{(-1)}),
 (1,\mathbf{0}) ) \in K^*(2),
\end{equation*}

\noindent such that $Q_1 \ _1\infty_{-2} \ Q_2$ exists. The unique
$F^{(1)}$ and $F^{(2)}$ satisfying (\ref{E:sewing_eqn1}) are as in
the previous example.  Similarly, define $\hat{f}_2$ as above.
%\begin{align*}
%F^{(1)}(w)=&w \\
%F^{(2)}(w)=& f_1^{-1}(w)= (a_0^{(1)})^{-w \frac{d}{dw}}
%e^{\sum_{j \in \Z} A_+j^{(1)} w^{j+1} \frac{d}{dw}}w.
%\end{align*}
Thus

\begin{equation} \label{E:ex5}
Q_1 \ _1\infty_{-2} \ Q_2 = ((a_0^{(1)}\hat{f}_2(\zeta))^{-1};
A^{(-1)},
%include to nontrivialize  B^{(1)} (c_0^{(1)}, C^{(1)}),
(c_0^{(-1)}, C^{(-1)}), (a_0^{(1)}, A^{(1)}) )
\end{equation}

\noindent
%include to make nontrivial
%where $c_0^{(i)}$, $C^{(i)}=\{C^{(i)}_j\}_{j \in \Z+}$, for $i=1,-1$
where $c_0^{(-1)}$, $C^{(-1)}=\{C^{(-1)}_j\}_{j \in \Z+}$ are
defined by

%include to make nontrivial
%\begin{multline*}
%e^{\sum_{j \in \Z_+} C_j^{(1)} w^{j+1} \frac{d}{dw}}
%(c_0^{(1)})^{w \frac{d}{dw}} \\
%=e^{\sum_{j \in \Z_+} A_j^{(1)} w^{j+1} \frac{d}{dw}}
%(a_0^{(1)}^{w \frac{d}{dw}}
%e^{\sum_{j \in \Z_+} B_j^{(1)} w^{j+1} \frac{d}{dw}}
%(b_0^{(1)})^{w \frac{d}{dw}},
%\end{multline*}
\begin{multline*}
e^{-\sum_{j \in \Z_+} C_j^{(-1)} w^{-j+1} \frac{d}{dw}}
(c_0^{(-1)})^{-w \frac{d}{dw}} \\
=e^{-\sum_{j \in \Z_+} (a_0^{(1)})^{-j} \Theta_j^{(2)} w^{-j+1}
\frac{d}{dw}} \left(\frac{e^{\Theta_0^{(2)}}}{a_0^{(1)}}
\right)^{-w \frac{d}{dw}} e^{-\sum_{j \in \Z_+} B_j^{(-1)}
w^{-j+1} \frac{d}{dw}} (b_0^{(-1)})^{-w \frac{d}{dw}}.
\end{multline*}

\noindent Note that $c_0^{(-1)}$ and $C^{(-1)}$ are uniquely
determined (and independent of Virasoro representation) by
\cite{BHL}.
\end{example}

We could also consider the sewing of $Q_1=(z^{-1};w^{-1},
w^{-1}-z, w) \in K^*(2)$ and $Q_2=(B^{(-1)}, (1,\mathbf{0}) ) \in
K^*(1)$ such that $Q_1 \ _1\infty_{-1} \ Q_2$ exists, and
similarly achieve another sewing identity involving
$\Theta_j^{(1)}$'s and $\hat{f}_1(x)=e^{\sum_{j \in \Z}
\mathcal{B}_j^{(-1)} x^{j+1} \frac{d}{dx}}x$:

\begin{multline} \label{E:firstsewingidentity}
e^{\sum_{k=-1}^{\infty} \left( \sum_{j \in \Z_+} { j+1 \choose
k+1} \mathcal{B}^{(-1)}_j z^{j-k}
\right)w^{-k+1} \frac{\partial}{\partial w}} \\
=e^{ \Theta^{(1)}_0 w\frac{\partial}{\partial w}} e^{\sum_{j \in
\Z_+} \Theta_j^{(1)} w^{-j+1} \frac{\partial}{\partial w}}
e^{(z-\hat{f}_2^{-1}(z))w^2 \frac{\partial}{\partial w}}.
\end{multline}

\noindent (called the first sewing identity).  This is done
explicitly in \cite{K} but directly follows the strategy we have
used above.

\subsection{The linear functionals $\mathcal{L}_I(z)$} \label{S:L_I}
In our work to understand representations of the Virasoro algebra
in a given vertex operator coalgebra, a particular family of
linear functionals will be of interest.  First, we define the
notion of a meromorphic function on any manifold $K^*(n)$, for $n
\in \N$.

\begin{definition}
A meromorphic function on $K^*(n)$ is a function that, when viewed
as a function of $z_{-i}$, for $i=1, \ldots, n-1$, of
$a_0^{(-i)}$, for $i=-1,1, \ldots, n-1$, and of $A_j^{(-i)}$, for
$i= -1,1, \ldots, n$ and $j \in \Z_+$, is a polynomial in these
variables divided by the product of powers of $a_0^{(-i)}$, for
$i=-1,1, \ldots, n-1$, of $z_{-i}$, for $i=1, \ldots, n-1$, and of
$z_{-i}-z_{-j}$, for $1 \leq i < j \leq n-1$ (cf. \cite{H} p. 67).
Explicitly, if $F$ is a meromorphic function on $K^*(n)$ and $Q
\in K^*(n)$ is as in (\ref{E:canon1}) or (\ref{E:canon2}), then

\begin{multline*}
F(Q) = \frac{\tilde{F}(z_{-n+1}, \ldots, z_{-1},A^{(-n)}, \ldots,
A^{(-1)}, A^{(1)}, a_0^{(-n+1)}, \ldots, a_0^{(-1)}, a_0^{(1)} ) }
{(a_0^{(1)})^{r_1}\prod_{i=1}^{n-1} (a_0^{(-i)})^{r_{-i}} \prod_{i
<j} (z_{-i}-z_{-j})^{s_{ij}} \prod_{i=1}^{n-1} z_{-i}^{t_{-i}}},
\end{multline*}

\noindent where $\tilde{F}$ is a polynomial, and $r_i, s_{ij}, t_i
\in \N$.
\end{definition}

In addition to the definition, we will also use the following
proposition which is Proposition 2.1.17 in \cite{H}.

\begin{proposition} \label{P:compose}
Let $g(x) \in R[[x,x^{-1}]]$ and $f(x) \in R[[x]]$ with

\begin{equation*}
f(x) =e^{ \sum_{k \in \N} \mathcal{A}_k x^{k+1} \frac{d}{dx}}x.
\end{equation*}

\noindent Then if $g \circ f$ is well-defined,

\begin{equation*}
g \circ f(x) = e^{\sum_{k \in \N} \mathcal{A}_k x^{k+1}
\frac{d}{dx}} g(x).
\end{equation*}
\end{proposition}

The linear functionals of interests will be denoted
$\mathcal{L}_I(z)$, $z \in \C^{\times}$, and will map from the
space of meromorphic functions on $K^*(1)$ to $\C$.  We define
$\mathcal{L}_I(z)$ by

\begin{equation*}
\mathcal{L}_I(z)F=\left. \left(\frac{d}{d \ep} F \left(
((1,(0,-\ep,0,0,\ldots)) ) _1\infty_{-1} (z^{-1} ;\mathbf{0},
(1,\mathbf{0}),(1,\mathbf{0})) \right) \right) \right|_{\ep=0}.
\end{equation*}

Since $\mathcal{L}_I(z)$ may be thought of as an element of the
tangent space of $K^*(1)$ at $(\mathbf{0}, (1,\mathbf{0}))$, it
may be rewritten

\begin{equation*}
\mathcal{L}_I(z)= \left. \sum_{k \in \Z_+} \gamma_{-k}
\frac{\partial}{\partial \mathcal{A}_k^{(-1)}}
\right|_{\mathcal{A}^{(-1)}= \mathbf{0}}+ \left. \gamma_0
\frac{\partial}{\partial \alpha_0^{(1)}}
\right|_{\alpha_0^{(1)}=1} +\left. \sum_{k \in \Z_+} \gamma_k
\frac{\partial}{\partial \mathcal{A}_k^{(1)}}
\right|_{\mathcal{A}^{(1)}= \mathbf{0}}
\end{equation*}

\noindent for some coefficients $\gamma_k$, with $k \in \Z$.  In
fact, we can completely describe these coefficients as the
following lemma shows.

\begin{lemma} \label{L:L_I}
\begin{equation*}
\mathcal{L}_I(z)= \sum_{k \in \Z_+} -z^{-k-2} \left.
\frac{\partial}{\partial \mathcal{A}_k^{(-1)}}
\right|_{\mathcal{A}^{(-1)}= \mathbf{0}} - \left. z^{-2}
\frac{\partial}{\partial \alpha_0^{(1)}}
\right|_{\alpha_0^{(1)}=1} +\sum_{k \in \Z_+} \left. -z^{k-2}
\frac{\partial}{\partial \mathcal{A}_k^{(1)}}
\right|_{\mathcal{A}^{(1)}= \mathbf{0}}.
\end{equation*}
\end{lemma}

\begin{proof}
The idea of the proof will be to sew $((1,(0,-\ep,0,0,...)) )$ and \\
$(z^{-1};\mathbf{0}, (1,\mathbf{0}),(1,\mathbf{0}))$ together then
examine what the derivative with respect to $\ep$ at $\ep =0$ is
for each of the local coordinate coefficients.  In some ways this
reflects the approach in \cite{H} to a similar linear functional
$\mathcal{L}_I(z)$, but here we are forced to consider nontrivial
scaling at the outgoing puncture and also the technique we use
avoids the unresolved branching issue in (3.2.10) and (3.2.11) of
\cite{H}.
%The following diagram will be helpful.

%\medskip
%DIAGRAM
%insert diagram here (:
%begin{figure}
%        \centering
%        \includegraphics[height=200pt]{figures/ex1.eps}
%        \caption{Puncture and Local Coordinate}
%        \label{fig_puncture}
%\end{figure}

%$$\epsfxsize=4.2in\epsfbox{figures/ex1.eps}$$

First, we denote the local coordinate map of
$((1,(0,-\ep,0,0,...)) )$ as $f_1(w)$ and the local coordinate
maps of $(z^{-1};\mathbf{0}, (1,\mathbf{0}),(1,\mathbf{0}) )$ as
$g_i(w)$ where the subscript indicates which puncture the local
coordinate map pertains to. Thus we get the sewing maps
$F^{(1)}(w)=w$ and
$F^{(2)}=f_1^{-1}\left(\frac{1}{g_{-1}(w)}\right)=
e^{zw^2\frac{d}{dw}}e^{\ep w^3\frac{d}{dw}}w$.  These maps move
the last outgoing puncture to $k=\left. e^{\ep x^3\frac{d}{dx}}x
\right|_{x=-z^{-1}}$.  We employ the transformation map
$T_1(w)=e^{k^{-1}w^2\frac{d}{dw}}w$.  (Notice this does not
guarantee that we will have a canonical representative of an
element in $K^*(1)$ since there may be a nontrivial linear scaling
of the local coordinate map at $\infty$.) After sewing and
applying $T_1$, the local coordinate map at 0 is $g_1 \circ
(F^{(2)})^{-1} \circ T_1^{-1}(w)$ and the local coordinate map at
$\infty$ is $g_{-2} \circ (F^{(2)})^{-1} \circ T_1^{-1}(w)$, both
of which depend on $\ep$ and will be trivial if $\ep=0$.

Now, we take the derivative of the puncture at 0 with respect to
$\ep$ and evaluate at $\ep=0$ (using Proposition \ref{P:compose}
twice).

\begin{equation*}
\frac{d}{d \ep} \left.\left(g_1 \circ (F^{(2)})^{-1} \circ
T_1^{-1}(w) \right) \right|_{\ep=0} =\frac{d}{d \ep}
\left.\left(e^{-k^{-1}w^2\frac{d}{dw}} e^{-\ep
w^3\frac{d}{dw}}e^{-zw^2\frac{d}{dw}}w \right) \right|_{\ep=0}
\end{equation*}
\begin{equation*}
=\frac{d}{d \ep} \left.\left(e^{-(-z^{-1}-\ep
z^{-3}-\ldots)^{-1}w^2\frac{d}{dw}} ((w^{-1}+z)^{-1}-\ep w
(w^{-1}+z)^{-2} + \ldots) \right) \right|_{\ep=0}
\end{equation*}
\begin{multline}\label{E:103}
=\frac{d}{d \ep} \left( ((w^{-1}-(z^{-1}+\ep
z^{-3}+\ldots)^{-1})+z)^{-1}
-\ep (w^{-1}- \right. \\
\left.\left. (z^{-1}+\ep z^{-3}+\ldots)^{-1})^{-1}
((w^{-1}-(z^{-1}+\ep z^{-3}+\ldots)^{-1})+z)^{-2} + \ldots \right)
\right|_{\ep=0}
\end{multline}

\noindent We only use ``$\ldots$" for terms involving higher
powers of $\ep$, which will disappear upon differentiating and
evaluating at $\ep=0$, so this notation is unambiguous.  Now
applying the derivative and evaluating, the right-hand side of
(\ref{E:103}) at $\ep=0$, we obtain

\begin{equation*}
-(w^{-1}-z+z)^{-2}z^2z^{-3}-(w^{-1}-z)^{-1} (w^{-1}-z+z)^{-2}
\end{equation*}
\begin{align*}
&=-w^{2}z^{-1}-w^2(w^{-1}-z)^{-1} \\
&=-w^2\frac{z^{-1}}{1-wz}
\end{align*}

We are dealing with local coordinates at 0 which means there will
be no negative powers of $w$. Thus

\begin{equation*}
\frac{d}{d \ep} \left.\left(g_1 \circ (F^{(2)})^{-1} \circ
T_1^{-1}(w) \right) \right|_{\ep=0} =-w^2 \frac{z}{1-wz}
\end{equation*}
\begin{equation} \label{E:110}
=-\sum_{k=0}^{\infty} z^{k-1}w^{k+2}.
\end{equation}

Using the same approach for the local coordinates at $\infty$, we
find that

\begin{equation*}
\frac{d}{d \ep} \left.\left(g_{-2} \circ (F^{(2)})^{-1} \circ
T_1^{-1}(w) \right) \right|_{\ep=0} =\frac{d}{d \ep}
\left.\left(\frac{1}{ (F^{(2)})^{-1} \circ T_1^{-1}(w)} \right)
\right|_{\ep=0}
\end{equation*}
\begin{equation*}
=\left. -((F^{(2)})^{-1} \circ T_1^{-1}(w))^{-2} \frac{d}{d \ep}
\left( (F^{(2)})^{-1} \circ T_1^{-1}(w) \right) \right|_{\ep=0}
\end{equation*}
\begin{equation*}
=-w^{-2}\left(-w^2\frac{z^{-1}}{1-wz}\right)
\end{equation*}
\begin{equation*}
=-\frac{w^{-1}z^{-2}}{1-w^{-1}z^{-1}}
\end{equation*}
\begin{equation}\label{E:111}
=-\sum_{k=0}^{\infty} z^{-k-2}w^{-k-1}.
\end{equation}

We see now that the local coordinates at $\infty$ do have an $a_0$
term that depends on $\ep$ since their derivative evaluated at
$\ep=0$ has a nontrivial coefficient of $w^{-1}$ in the power
series expansion, and we also see that since the coefficient of
$w$ in the power series expansion is 0, the local coordinates at 0
have a trivial $a_0$ term. We want to find out what the
coefficients are in front of the canonical tangent vectors,
though, so let $h_1(w)=e^{\sum_{k \in \Z_+} A^{(1)}_k(\ep)
w^{k+1}\frac{d}{dw}} (a_0^{(1)}(\ep))^{w\frac{d}{dw}}w$ and
$h_{-1}(w)=e^{\sum_{k \in \Z_+} -A^{(-1)}_k(\ep)
w^{-k+1}\frac{d}{dw}}\frac{1}{w}$ be the canonical representatives
of the local coordinates obtained from the sewing (where
$A^{(i)}_k(\ep), a_0^{(1)}(\ep)$ are functions depending on $\ep$
for $i=-1,1$ and $k \in \Z_+$). By applying the transformation map
$T_2(w)=(a_0^{(1)}(\ep))^{w\frac{d}{dw}}w$, we again obtain the
unique local coordinates with no nontrivial $a_0$ at 0 and a
possibly nontrivial $a_0$ term at $\infty$.  Taking the derivative
of these local coordinates with respect to $\ep$ at $\ep=0$ yields

\begin{equation}\label{E:112}
\frac{d}{d \ep} \left. \left(h_1 \circ T_2^{-1}(w) \right)
\right|_{\ep=0}
=\sum_{k \in \Z_+} \left. \left( \frac{d}{d \ep}
A^{(1)}_k(\ep) \right) \right|_{\ep=0} w^{k+1}
\end{equation}

\noindent for the puncture at 0 and

\begin{equation*}
\frac{d}{d \ep} \left.\left(h_{-1} \circ T_2^{-1}(w) \right)
\right|_{\ep=0} =\frac{d}{d \ep} \left.\left(
(a_0^{(1)}(\ep))^{-w\frac{d}{dw}} e^{-\sum_{k \in \Z_+}
A^{(-1)}_k(\ep) w^{-k+1}\frac{d}{dw}}\frac{1}{w} \right)
\right|_{\ep=0}
\end{equation*}
\begin{equation}\label{E:113}
=\left.\left( \frac{d}{d \ep}  (a_0^{(1)}(\ep))\right)
\right|_{\ep=0}w^{-1} + \sum_{k \in \Z_+} \left.\left( \frac{d}{d
\ep} A^{(-1)}_k(\ep) \right) \right|_{\ep=0}w^{-k-1}
\end{equation}

\noindent for the puncture at $\infty$.

By comparing the coefficients of (\ref{E:110}) and (\ref{E:112})
as well as (\ref{E:111}) and (\ref{E:113}), we verify the claim of
the lemma.
\end{proof}

\section{Algebraic preliminaries}
It is a standard notion in the theory of operads to study the
formal algebraic structure that a geometric structure induces on
vector spaces (cf. \cite{M}, \cite{MSS} and \cite{HL}). In order
to study the kind of structure that $K^*$ induces on vector
spaces, it is necessary to define and investigate some algebraic
preliminaries.  This section will cover the formal calculus we
will need. For a more thorough exposition on formal calculus
itself, see
 \cite{FLM}, \cite{FHL} or \cite{LL}.

\subsection{The $\delta$-function}

We will use the ``formal $\delta$-function", $\delta(x)= \sum_{n
\in \Z} x^n$, which is discussed in, for instance, \cite{FHL}.
Note that the $\delta$-function applied to $\frac{x_1-x_2}{x_0}$,
where $x_0$, $x_1$ and $x_2$ are commuting formal variables, is a
formal power series in $x_2$ (i.e., negative powers of $(x_1-x_2)$
are expanded in nonnegative integral powers of $x_2$).  In
general, in the formal calculus of VOAs the sum or difference of
two formal variables, $(x_1 \pm x_2)$, is understood to be
expanded in nonnegative integral powers of the second, $x_2$ (cf.
\cite{FLM}, \cite{FHL}).

The following three properties of $\delta$-functions will be
relevant:

\noindent First, given a formal Laurent series $X(x_1,x_2) \in
\text{Hom} (V,W)[[x_1,x_1^{-1},x_2,x_2^{-1}]]$ with coefficients
which are homomorphisms from a vector space $V$ to a vector space
$W$, if $\displaystyle \lim_{x_1 \to x_2} X(x_1,x_2)$ exists (i.e.
when $X(x_1,x_2)$ is applied to any element of $V$, setting
$x_1=x_2$ leads to only finite sums in $W$) we have

\begin{equation}\label{E:delta1}
\delta \left( \frac{x_1}{x_2} \right) X(x_1,x_2) = \delta
\left(\frac{x_1}{x_2} \right) X(x_2,x_2).
\end{equation}

%(should I include this ?)
%Another useful form of this identity is:
%\begin{equation*}
%\delta \left( \frac{x_1-x_0}{x_2} \right) X(x_0,x_1,x_2)
%= \delta \left(\frac{x_1-x_0}{x_2} \right)\iota_{1 \ 0} X(x_0,x_1,x_1-x_0).
%\end{equation*}

\noindent Second, we use the fact that

\begin{equation}\label{E:delta2}
x_1^{-1}\delta \left(\frac{x_2+x_0}{x_1}  \right) =x_2^{-1}\delta
\left(\frac{x_1-x_0}{x_2} \right)
\end{equation}

\noindent which is proved by direct expansion.  The third fact may
be observed by direct expansion and comparing of coefficients:

\begin{equation}\label{E:delta3}
x_0^{-1}\delta \left(\frac{x_1-x_2}{x_0}  \right) -x_0^{-1}\delta
\left(\frac{x_2-x_1}{-x_0} \right) =x_2^{-1}\delta
\left(\frac{x_1-x_0}{x_2} \right).
\end{equation}

\subsection{Linear algebra on $\Z$-graded vector spaces with finite-dimensional
homogeneous subspaces}

Let $V = \coprod_{n \in \Z} V_{(n)}$ be a $\Z$-graded vector space
over $\C$ such that $\dim V_{(n)} < \infty$, for each $n \in \Z$.
We denote the graded dual space of $V$ by

\begin{equation*}
V'= \coprod_{n \in \Z} V^*_{(n)},
\end{equation*}

\noindent
the algebraic closure of $V$ by

\begin{equation*}
\overline{V}= \prod_{n \in \Z} V_{(n)}=(V')^*,
\end{equation*}

\noindent
and the natural pairing of $V'$ with $\overline{V}$ by $\langle \cdot, \cdot \rangle$.
The $n$-th tensor product of $V$, denoted $V^{\otimes n}$, is still a
$\Z$-graded vector space (where $v \in V_{(k_1)} \otimes \ldots \otimes V_{(k_n)}$ has
weight $k_1 + \ldots + k_n$) with finite-dimensional homogeneous subspaces.  Thus
$(V^{\otimes n})'$, $\overline{V^{\otimes n}}$ and  $\langle \cdot, \cdot \rangle:
(V^{\otimes n})' \times \overline{V^{\otimes n}} \to \C$
are defined as above.

We denote a homogeneous basis of $V$ by

\begin{equation*}
\{ e^{(k)}_{\ell^{(k)}} | k \in \Z, \ \ell^{(k)}=1, \ldots, \dim V_{(k)} \}
\end{equation*}

\noindent
and its corresponding dual basis of $V'$ by

\begin{equation*}
\{ (e^{(k)}_{\ell^{(k)}})^* | k \in \Z, \ \ell^{(k)}=1, \ldots, \dim V_{(k)} \}.
\end{equation*}

\noindent
We will use the notation

\begin{equation*}
\mathcal{H}_V(m,n)=\text{Hom}(V^{\otimes m}, \overline{V^{\otimes
n}})
\end{equation*}

\noindent for $m,n \in \N$.  For $m \in \N$, $n \in \Z_+$ and any
integer $0 < i \leq n$, we define the \emph{t-contraction}

\begin{align*}
( \cdot \ _1*_{-i} \ \cdot)_t: \mathcal{H}_V(1,m) \times
\mathcal{H}_V(1,n)
&\to \text{Hom}(V, \overline{V^{\otimes m+n-1}}[[t,t^{-1}]])\\
(f,g) &\mapsto (f \ _1*_{-i} \ g)_t
\end{align*}

\noindent
where

\begin{multline*}
(f \ _1*_{-i} \ g)_t v = \sum_{k \in \Z} \sum_{\ell^{(k)}=1}^{\dim V_{(k)}}
\left(\underbrace{Id_V \otimes \cdots \otimes Id_V}_{i-1} \otimes
f (e^{(k)}_{\ell^{(k)}}) \cdot
 (e^{(k)}_{\ell^{(k)}})^*\otimes \right. \\
\left. \underbrace{Id_V \otimes \cdots \otimes Id_V}_{n-i}\right) g(v) t^k.
\end{multline*}

\noindent
If for all $v' \in (V^{\otimes n+m-1})'$, $v\in V$ the formal Laurent series

\begin{equation*}
\langle v',(f \ _1*_{-i} \ g)_t v \rangle
\end{equation*}

\noindent is absolutely convergent when $t=1$, then $(f \ _1*_{-i}
\ g)_1$ is well-defined as an element of $\mathcal{H}_V(1,m+n-1)$
and we define the \emph {contraction} of $f$ and $g$ by,

\begin{equation*}
f \ _1*_{-i} \ g= (f \ _1*_{-i} \ g)_1.
\end{equation*}

The following associativity of $t$-contractions follows from the definition.

\begin{proposition} \label{P:H^*_assoc}
Let $\ell,m \in \N$ and $n \in \Z_+$ such that $m+n > 1$.  Choose
$f_1 \in \mathcal{H}_V(1,\ell), \ f_2 \in \mathcal{H}_V(1,m), \
f_3 \in \mathcal{H}_V(1,n)$ and integers $1 \leq i \leq m+n-1$ and
$1 \leq j \leq n$ .  One of the following 3 holds:

\noindent
(1) $i < j$ and as a formal series in $t_1$ and $t_2$

\begin{equation*}
(f_1 \ _1*_{-i} \ (f_2 \ _1*_{-j} \ f_3)_{t_1})_{t_2}
=(f_2 \ _1*_{-j-\ell+1} \ (f_1 \ _1*_{-i} \ f_3)_{t_2})_{t_1};
\end{equation*}

\noindent
(2) $j \leq i < j+m$ and as a formal series in $t_1$ and $t_2$

\begin{equation*}
(f_1 \ _1*_{-i} \ (f_2 \ _1*_{-j} \ f_3)_{t_1})_{t_2}
=((f_1 \ _1*_{-i+j-1} \ f_2)_{t_2} \ _1*_{-j} \ f_3)_{t_1};
\end{equation*}

\noindent
(3) $i \geq j+m$ and as a formal series in $t_1$ and $t_2$

\begin{equation*}
(f_1 \ _1*_{-i} \ (f_2 \ _1*_{-j} \ f_3)_{t_1})_{t_2}
=(f_2 \ _1*_{-j} \ (f_1 \ _1*_{-i+m-1} \ f_3)_{t_2})_{t_1}.
\end{equation*}
\end{proposition}

\noindent
Since this proposition shows equality as formal power series, it also
implies absolute convergence of both sides and equality given the absolute
convergence of either side.

Moving from associativity to permutations, the symmetric group on $n$
letters acts naturally on
$\overline{V^{\otimes n}}$ from the left, i.e. it is determined by

\begin{equation*}
\sigma(v_1 \otimes \ldots \otimes v_n)= v_{\sigma^{-1}(1)} \otimes v_{\sigma^{-1}(n)}
\end{equation*}

\noindent for all $\sigma \in S_n$ and $v_1, \ldots, v_n \in V$.
This induces a left action on $\mathcal{H}_V(1,n)$ given by

\begin{equation*}
\sigma(f)(v)=\sigma(f(v)).
\end{equation*}

\noindent Transpositions play a fundamental role in the actions of
the symmetric group so we will make heavy use of the transposition
map

\begin{align*}
T: V \otimes V &\to V \otimes V \\
 v \otimes w &\mapsto w \otimes v.
\end{align*}

\noindent Hence, from the definition of $t$-contraction and the
action of the symmetric group we see the following.

\begin{proposition} \label{P:H^*perm}
Let $f_1, f_2 \in \mathcal{H}_V(1,2)$. Then $f_1 \ _1\infty_{-1} \
f_2$ exists if and only if $(Id_V \otimes T)(T \otimes Id_V)(f_{1}
\ _1\infty_{-2} \ (T \ f_2))$ also exists.  If this is the case,

\begin{equation*}
f_{1} \ _1*_{-1} \ f_2= (Id_V \otimes T)(T \otimes Id_V)(f_{1} \ _1*_{-2} \ (T \ f_2)).
\end{equation*}
\end{proposition}

Proposition \ref{P:H^*perm} implies that all symmetric groups act functorially
with respect to contraction since every symmetric group is generated
by transpositions.

\begin{remark}
In his book \cite{H}, Huang examines a similar $t$-contraction and
contraction on the set of linear maps $\mathcal{H}_V(m,1)$ for $m
\in \Z$. The vector spaces $\mathcal{H}_V(m,1)$ and
$\mathcal{H}_V(1,m)$ are naturally isomorphic, and it can even be
shown that the $t$-contraction and contraction in
$\mathcal{H}_V(m,1)$ correspond to the $t$-contraction and
contraction in $\mathcal{H}_V(1,m)$ under this isomorphism.
Geometric vertex operator algebras arise from considering maps
$\nu_m:K(m) \to H(m,1)$ such that the diagram

$$
\begin{CD}
K(m) \times K(n)     @> \nu_m  \times \nu_n>> H(m,1) \times H(n,1)\\
@V_i\infty_{-1} VV   @V _i*_{-1}VV\\
K(m+n-1)             @> \nu_{m+n-1}>> H(m+n-1,1)
\end{CD}
$$

%\noindent and

%$$
%\begin{CD}
%K(n)  @> \nu_n>> H(n,1)\\
%@V \sigma VV  @V \sigma VV\\
%K(n) @> \nu_{n}>> H(n,1).
%\end{CD}
%$$

\noindent commutes (up to scalars, when the map $_i\infty_{-1}$ is
defined).
%, and for $\sigma \in S_n$).
Among other conditions for GVOAs, the grading and meromorphicity
axioms place additional requirements on the image of $\langle
v',\nu_m(Q)v \rangle$ that depend on the choice of elements $Q\in
K(m)$ and $v \in V^{\otimes m}$, but not on $v' \in V'$.  In some
sense we could use the isomorphisms $K(m) \cong K^*(m)$ (which are
not canonical - recall Remark \ref{R:operad_iso}) and $H(m,1)
\cong H(1,m)$ to achieve the diagram

$$
\begin{CD}
K^*(m) \times K^*(n)   @> \mu_m  \times \nu_n >> H(1,m) \times H(1,n)\\
@V_{1}\infty_{-i} VV   @V _{1}*_{-i}VV \\
K^*(m+n-1)             @> \mu_{m+n-1}>> H(1,m+n-1).
\end{CD}
$$

\noindent  This diagram underlies geometric vertex operator
coalgebras.  However, in order to maintain consistency with the
current worldsheet model, GVOCs must still place additional
requirements on the image of $\langle v',\mu_m(Q)v \rangle$ that
depend on the choice of elements $Q\in K^*(m)$ and $v \in V$, but
not on $v' \in (V^{\otimes m})'$.  These requirements would be
dualized under the above isomorphisms and hence would be
incorrect.  The  approach we provide in Section \ref{S:iso1}
yields a constructive isomorphism with a constructive inverse
(that is natural with respect to the formal definition of $K^*$).
The above approach is both nonconstructive and not canonical so
that, even if the discrepancies with the grading and
meromorphicity axioms were resolved, our approach produces a more
useful outcome.
\end{remark}

\subsection{More on sewing identities in representations of the Virasoro algebra}
Now that we have developed another object on which the Virasoro
algebra may act, it is important to highlight the independence of
the choice of the $\Theta$ sequences from the representation in
which they are selected. Without specifying the specific action of
the Virasoro algebra on a complex vector space $V$, we can claim
the sewing identities we discussed in Section \ref{S:example}
((\ref{E:secondsewingidentity}) and
(\ref{E:firstsewingidentity})). They correspond to Propositions
4.3.9 and 4.3.10 from \cite{H} with minor modifications.

\begin{proposition} \label{P:sew_id3b}
Let $\Theta^{(2)}=\Theta^{(2)}
(\mathcal{A},\alpha_0^{(1)},\zeta)$,
$\Theta_0^{(2)}=\Theta_0^{(2)}(\mathcal{A},\alpha_0^{(1)},\zeta)$
be chosen as in Proposition \ref{P:sew_id3} to depend on
$\mathcal{A}$, $\alpha_0^{(1)}$, and $\zeta$.

Then in the algebra $(\text{End } V)[[t]][y,y^{-1}]
[\alpha_0^{(1)}, \alpha_0^{(1)}] [[\mathcal{A}^{(1)}]]$ we have

\begin{multline*}
e^{\sum_{k=-1}^{\infty} \left( \sum_{j \in \Z_+} (\alpha_0^{(1)}
)^{-j} \mathcal{A}^{(1)}_j
{ -j+1 \choose k+1}\zeta^{-j-k} \right)L(-k)}\\
=e^{- \Theta^{(2)}_0 L(0)} e^{-\sum_{j \in \Z_+} \Theta_j^{(2)}
L(-j)} e^{(-\zeta + \hat{f}_2(\zeta))L(1)} .
\end{multline*}
\end{proposition}

\begin{proposition} \label{P:sew_id2b}
Similarly, let $\Theta^{(1)}=\Theta^{(1)} (\mathcal{B},z)$,
$\Theta_0^{(1)}=\Theta_0^{(1)} (\mathcal{B},z)$ be as in Equation
\ref{E:firstsewingidentity}. Then in the algebra $(\text{End }
V)[[t]][y,y^{-1}][[\mathcal{B}^{(-1)}]]$ we have

\begin{multline*}
e^{-\sum_{k=-1}^{\infty} \left( \sum_{j \in \Z_+} { j+1 \choose
k+1} \mathcal{B}^{(-1)}_j z^{j-k}
\right) L(-k)} \\
= e^{-\Theta^{(1)}_0 L(0)} e^{-\sum_{j \in \Z_+} \Theta_j^{(1)}
L(-j)} e^{(-z+\hat{f}_1^{-1}(z)) L(1)}.
\end{multline*}
\end{proposition}

\section{The notion of vertex operator coalgebra and its geometric interpretation}
We now have the background to define the notion of (algebraic)
vertex operator coalgebra axiomatically as well as the notion of
geometric vertex operator coalgebra.

\subsection{The notion of geometric vertex operator coalgebra}
We begin by defining the primary geometric motivation for vertex
operator coalgebras. The definition of a geometric vertex operator
coalgebra uses the same kind of conformal structure that underlies
the definition of a geometric vertex operator algebra first given
in \cite{H2} and \cite{H} (in order that they might eventually be
combined). One may think of a geometric vertex operator algebra as
a meromorphic morphism (or algebra) associated to a $\C$-extension
of the partial operad $K$. (See Section 6 of \cite{HL}.)
Similarly, a geometric vertex operator coalgebra may be
interpreted as a meromorphic morphism (or coalgebra) of a
$\C$-extension of the partial operad $K^*$ to the partial
pseudo-operad $\{ \mathcal{H}_V (1,n) \}_{n \in \N}$. We will not
focus on this interpretation here.

\begin{definition}\label{GVOC}
    A \emph{geometric vertex operator coalgebra (over $\C$) of rank $d \in \C$} is
a $\Z$-graded vector space over $\C$

\begin{equation*}
V = \coprod_{k \in \Z} V_{(k)}
\end{equation*}

\noindent such that \emph{dim}$V_{(k)} < \infty$ for $k \in \Z$
and $V_{(k)} = 0$ for $k$ sufficiently small, together with a
linear map for each $n \in \N$

\begin{equation*}
\mu_n : K^*(n) \mapsto \mathcal{H}_V (1,n)
\end{equation*}

\noindent
satisfying the following axioms:

1. Grading: Let $k \in \Z$, $v \in V_{(k)}$, $v' \in V'$, and $a \in \C^\times$. Then

\begin{equation*}
\langle v', \mu_1(\mathbf{0},(a,\mathbf{0}))v \rangle= a^{-k}
\langle v',v \rangle.
\end{equation*}

2. Meromorphicity: For any $n \in \N$, $v' \in (V^{\otimes n})'$, and $v \in V$ the function

\begin{equation*}
Q \mapsto \langle v', \mu_n(Q)v \rangle
\end{equation*}

\noindent
from $K^*(n)$ to $\C$
is a meromorphic function on $K^*(n)$.  Further, for any $v \in V$ there exists
$N(v) \in \Z_+$ such that for all $v' \in (V^{\otimes n})'$ the degree of
$z_{-i}$ in $\langle v', \mu_n(Q)v \rangle$ is less than $N(v)$, for $i=1, \ldots, n-1$.

3. Permutation: Let $\sigma \in S_n$.  Then for any $Q \in K^*(n)$

\begin{equation*}
\sigma(\mu_n(Q)) = \mu_n(\sigma(Q)).
\end{equation*}

4. Sewing: Given

\begin{multline*}
Q_1=(z_{-m+1}^{-1}, \ldots, z_{-1}^{-1}; A^{(-m)}, (a_0^{(-m+1)},
A^{(-m+1)}), \ldots,
(a_0^{(-1)},A^{(-1)}), \\
(a_0^{(1)},A^{(1})) \in K^*(m),
\end{multline*}
\begin{multline*}
Q_2=(\zeta_{-n+1}^{-1}, \ldots, \zeta_{-1}^{-1}; B^{(-n)},
(b_0^{(-n+1)},B^{(-n+1)}), \ldots,
(b_0^{(-1)},B^{(-1)}), \\
(b_0^{(1)},B^{(1})) \in K^*(n),
\end{multline*}

\noindent
if the $i$-th outgoing puncture $(1 \leq i \leq n)$ of $Q_2$ can be sewn
to the incoming puncture of $Q_1$, then for any $v' \in (V^{\otimes m+n-1})',
v \in V$

\begin{equation*}
\langle v',(\mu_m(Q_1) \ _1*_{-i} \ \mu_n(Q_2)_t(v) \rangle
\end{equation*}

\noindent
is absolutely convergent when $t = 1$, and

\begin{equation*}
 \mu_{m+n-1} (Q_1 \ _1\infty_{-i} \ Q_2) =
(\mu_n(Q_1) _1*_{-i} \mu_m(Q_2) e^{-\Gamma(A^{(1)}, B^{(-i)},
a_0^{(1)}b_0^{(-i)})d}
\end{equation*}

\noindent where $d$ is the rank and $\Gamma(A^{(1)}, B^{(-i)},
a_0^{(1)}b_0^{(-i)})$ is as in (\ref{E:4162}).
\end{definition}

We denote the geometric vertex operator coalgebra just defined by
$(V, \mu = \{\mu_n\}_{n \in \N})$ or, when there is no ambiguity, simply $V$.

\begin{remark}
The second half of the meromorphicity axiom may be interpreted as
saying that $\langle v', \mu_n(Q)v \rangle$, when viewed as a
function in $z_{-1}^{-1}, \ldots, z_{-n+1}^{-1}$, has poles at
$z_{-i}^{-1}=0$ (for $i=1, \ldots, n-1$) and
$z_{-i}^{-1}=z_{-j}^{-1}$ (for $i \neq j$) such that the order of
each of these poles is bounded by $N(v)$ independent of $v' \in
(V^{\otimes n})'$. This fact comes directly from the
meromorphicity axiom for poles at $z_i^{-1}=0$ but for poles at
$z_i^{-1}=z_j^{-1}$ this observation involves the sewing axiom.
Incidently, the meromorphicity axiom of a geometric vertex
operator algebra includes a small bit of redundancy in that the
order of the poles at $z_i=z_j$ (for $i \neq j$) is bounded as a
result of the bounding of poles at $z_i=0$ (for $i=1, \ldots,
n-1$) and the sewing axiom.
\end{remark}

\subsection{The notion of vertex operator coalgebra}
The following description of a vertex operator coalgebra is the
central definition of this paper. It is in a sense the culmination
of our effort to algebraically understand the structure induced on
vector spaces by moduli spaces of spheres with outgoing tubes.  In
another sense, however, it is a starting point for algebraic
study.  A few algebraic properties will be described in this
section but these barely scratch the surface in reference to
questions about the structure of vertex operator coalgebras and
representations over them.  See \cite{K2} and \cite{K3} for
additional exploration into the algebraic properties of VOCs.

\begin{definition}\label{D:voc}
A \emph{vertex operator coalgebra (over $\C$) of rank $d \in \C$} is
a $\Z$-graded vector space over $\C$

\begin{equation*}
V = \coprod_{k \in \Z} V_{(k)}
\end{equation*}

\noindent
such that $\dim V_{(k)} < \infty$ for $k \in \Z$ and $V_{(k)} = 0$ for
$k$ sufficiently small,
together with linear maps

\begin{align*}
\co (x): V &\mapsto (V \otimes V)[[x,x^{-1}]] \\
v &\mapsto \co(x)v = \sum_{k\in \Z} \Delta_k(v) x^{-k-1},
\end{align*}

\begin{equation*}
 c : V \mapsto \C,
\end{equation*}

\begin{equation*}
\rho : V \mapsto \C,
\end{equation*}

\noindent
called the \emph{coproduct}, the \emph{covacuum map} and the \emph{co-Virasoro map}, respectively,
satisfying the following 7 axioms:

1. Left Counit: For all $v \in V$

\begin{equation}\label{E:counit}
(c \otimes Id_V) \co(x)v=v
\end{equation}

2. Cocreation: For all $v \in V$

\begin{equation} \label{E:cocreat1}
(Id_V \otimes c) \co(x)v \in V[[x]] \ \ \text{and}
\end{equation}
\begin{equation} \label{E:cocreat2}
\lim_{x \to 0} (Id_V \otimes c) \co(x)v=v.
\end{equation}

3. Truncation: Given $v \in V$, then $\Delta_k(v) = 0$ for $k$ sufficiently small.

4. Jacobi Identity:

\begin{multline} \label{E:Jac}
x_0^{-1}\delta \left(\frac{x_1-x_2}{x_0} \right)
(Id_V \otimes \co(x_2))  \co(x_1)
-x_0^{-1}\delta \left(\frac{x_2-x_1}{-x_0} \right)
(T \otimes Id_V)\\
 (Id_V \otimes \co(x_1))  \co(x_2)
=x_2^{-1}\delta \left(\frac{x_1-x_0}{x_2} \right)
(\co(x_0) \otimes Id_V)  \co(x_2).
\end{multline}

5. Virasoro Algebra:
The Virasoro algebra bracket,

\begin{equation*}
[L(j),L(k)]=(j-k)L(j+k)+\frac{1}{12}(j^3-j)\delta_{j,-k}d,
\end{equation*}

\noindent
holds for $j, k \in \Z$, where

\begin{equation} \label{E:L_def}
(\rho \otimes Id_V) \co(x) = \sum_{k \in \Z} L(k) x^{k-2}.
\end{equation}

6. Grading:  For each $k \in \Z$ and $v \in V_{(k)}$

\begin{equation} \label{E:VOCgrading}
L(0)v= kv.
\end{equation}

7. $L(1)$-Derivative:

\begin{equation} \label{E:L(1)deriv}
\frac{d}{dx} \co(x)=
(L(1) \otimes Id_V) \co(x).
\end{equation}

\end{definition}

\noindent We denote this vertex operator coalgebra by $(V, \co, c,
\rho)$, or sometimes just $V$.

Note that $x$, $x_0$, $x_1$ and $x_2$ are formal commuting variables and
$\co$ is linear so that, for example,
$(Id_V \otimes \co(x_1))$ acting on the coefficients of $\co(x_2)v \in
(V \otimes V)[[x_2,x_2^{-1}]]$ is well defined. Notice also, that when each
expression is applied to any element of $V$, the coefficient of each
monomial in the formal variables is a finite sum.

\begin{remark} \label{R:subtle}
Note that while the definition of vertex operator coalgebra is
stated in such a way as to maximize its resemblance to that of a
vertex operator algebra (as presented in \cite{FLM}, \cite{FHL},
and \cite{LL}) there are several important differences. First,
whereas the creation and truncation axioms of VOAs both bound the
power of the formal variable from below, in VOCs cocreation allows
only non-negative powers of $x$ while truncation allows only
finitely many positive powers of $x$ but infinitely many negative
powers.  Among other implications, this means that cocreation
actually generates polynomials in $x$. Second, note that the
representation of the Virasoro algebra has been inverted, then
shifted by $x^{-4}$. Finally, while the Jacobi identity
description highlights similarities between VOAs and VOCs,
examining the corresponding ``weak commutativity" properties
reveals substantial differences (cf. \cite{FHL} and \cite{K3}).
\end{remark}

\subsection{Properties of VOCs} \label{S:properties}
There are a number of interesting consequences of the above
axioms, most quite analogous to VOA properties.  Properties are
listed here without proof but the reader is encouraged to see
\cite{K} or \cite{K3} for more adequate justifications.

\begin{equation} \label{E:shift}
\left( e^{x_0 L(1)} \otimes Id_V \right) \co(x)= \co(x+x_0).
\end{equation}

\noindent
The \emph{$L(1)$-commutation formula}:

\begin{equation} \label{E:L(1)com}
\co(x_2) L(1) = (Id_V \otimes L(1)) \co(x_2) + (L(1) \otimes Id_V)  \co(x_2)
\end{equation}

\noindent
The \emph{$L(0)$-commutation formula}:

\begin{equation} \label{E:L(0)com}
\co(x_2) L(0) = (Id_V \otimes L(0)) \co(x_2) + (L(0) \otimes Id_V)  \co(x_2)
+ x_2(L(1) \otimes Id_V)\co(x_2)
\end{equation}

\noindent
The \emph{$L(-1)$-commutation formula}:

\begin{equation} \label{E:L(-1)com}
\co(x_2) L(-1) = (Id_V \otimes L(-1))  \co(x_2) +
\left( (x_2^2 L(1) +2 x_2 L(0) +L(-1)) \otimes Id_V \right) \co(x_2)
\end{equation}

\begin{multline} \label{E:rho_on_jac}
\co(x_2) \sum_{k \in \Z} L(k) x_1^{k-2}
- (Id_V \otimes \sum_{k \in \Z} L(k) x_1^{k-2})  \co(x_2)\\
=\mathrm{Res}_{x_0} x_2^{-1}\delta \left(\frac{x_1-x_0}{x_2}
\right) (\sum_{k \in \Z} L(k) x_0^{k-2} \otimes Id_V)  \co(x_2)
\end{multline}

\begin{equation}\label{E:distrib}
\left( e^{x_0 L(1)} \otimes Id_V \right) \co(x)
= \left( Id_V \otimes e^{-x_0 L(1)} \right) \co(x) e^{x_0 L(1)}
\end{equation}

\begin{equation}\label{E:cocreate2}
e^{x_0 L(1)} = (Id_V \otimes c)\co(x_0)
\end{equation}

\noindent
\emph{Anti-symmetry}:

\begin{equation}\label{E:antisym}
T\co(x)=  \co(-x) e^{x L(1)}
\end{equation}

\begin{equation} \label{E:cL(j)}
c L(j) = 0 \in \text{Hom}(V,\C) \text{    for } j \leq 1
\end{equation}

\begin{equation} \label{E:cL(2)}
c L(2) = \rho
\end{equation}

\begin{equation} \label{E:rhoL(0)}
\rho L(0) = 2\rho
\end{equation}

Several fundamental facts about VOCs necessitate an investigation of the weights of
vectors and operators.  Among other things, weights help us to understand the elements
$\Delta_n(v) \in V \otimes V$.  When $v \in V_{(k)}$, we say that $v$ is homogeneous
of \emph{weight} $k$, and the grading axiom says
$L(0)v=(\text{wt} \ v) v$.  Given an element $w\in V \otimes V$,
if $\left((L(0) \otimes Id_V) + (Id_V \otimes L(0))\right)w
=aw$ for some $a \in \C$, then we say that the weight of
$w$ is $a$.  This
agrees with the grading on $V \otimes V$ since for $v \in V_{(k)}$ and
$w \in V_{(\ell)}$, we have $v \otimes w \in (V \otimes V)_{(k+\ell)}$ and
$\text{wt } (v \otimes w) = k + \ell$.  With this motivation, for $v \in V_{(k)}$ and
$w \in V_{(\ell)}$, we define $x^{Id_V \otimes L(0)}(v \otimes w) = x^{\ell} (v \otimes w)$ and
$x^{L(0) \otimes Id_V}(v \otimes w) = x^k (v \otimes w)$.
For homogeneous vectors, we see that

\begin{equation} \label{E:weight}
\text{wt } \Delta_j (v)
=\text{wt } (v) +j+1.
\end{equation}

\noindent
Considering $v \in V_{(k)}$, it is now clear via (\ref{E:weight}) that for each $\Delta_j(v)$,

\begin{equation*}
\Delta_j (v) = \sum_{\ell=1}^{r_j} u_{\ell} \otimes w_{\ell}
\end{equation*}

\noindent
where each $u_{\ell} \in V_{(i_j)}$, each $w_{\ell} \in V_{( \text{wt } \Delta_j (v) - i_j)}$,
$r_j \in \Z_+$, and each $i_j \in \Z$.  We also have:

\begin{equation}
x^{-Id_V \otimes L(0)} \co(x) x_0^{L(0)}= x_0^{L(0) \otimes Id_V} \co(x_0 x)
\end{equation}

\begin{equation} \label{E:expL(0)}
(Id_V \otimes a^{-L(0)}) \co(x) a^{L(0)}= (a^{L(0)} \otimes Id_V) \co(a x) \text{ for }a \in \C^{\times}
\end{equation}

\begin{equation} \label{E:L(0)L(k)}
a^{L(0)}e^{bL(k)}=e^{a^{-k}bL(k)}a^{L(0)} \text{ for } a \in \C^{\times}, \ b \in \C, \ k \in \Z.
\end{equation}

\noindent
Given $v \in V_{(k)}$, the fact that

\begin{equation}
L(0) L(j)v= (k-j) L(j)v
\end{equation}

\noindent
for any $j,k \in \Z$, tells us that each operator $L(j)$
raises the weight of homogeneous vectors by $-j$ and is said to have
\emph{weight} $-j$, as is the case with VOAs.

The Jacobi identity may be replaced in the definition of VOC with properties called
\emph{right and left rationality, commutativity}, and \emph{associativity} as the following
two propositions show.  Let
$\iota_{1 \ 2}$ map rational functions in $x_1$ and $x_2$ to their expansion in only finitely
many negative powers of $x_2$.  Define $\iota_{2 \ 1}$ and $\iota_{2 \ 0}$ similarly (cf.
Section 3.1 of \cite{FHL}).

\begin{proposition} \label{P:RR,LR,com,ass}
Let $(V, \co, c, \rho)$ be a VOC.  For any $v' \in (V^{\otimes 3}
)'$ and $v \in V$ we have the following properties.

(i) Right rationality:

\begin{equation*}
\langle v', (Id_V \otimes \co(x_2)) \co(x_1)v \rangle =
\iota_{1 \ 2} \frac{g(x_1,x_2)}{x_1^r x_2^s (x_1-x_2)^t}
\end{equation*}

\noindent
for some $g(x_1,x_2) \in \C[x_1,x_2]$ and $r,s,t \in \Z$, where $g(x_1,x_2)$ is unique
up to choice of $r$, $s$ and $t$.

(ii) Left rationality:

\begin{equation*}
\langle v', (\co(x_0) \otimes Id_V) \co(x_2)v \rangle =
\iota_{2 \ 0} \frac{h(x_0,x_2)}{x_0^r x_2^s (x_0+x_2)^t}
\end{equation*}

\noindent
for some $h(x_0,x_2) \in \C[x_0,x_2]$ and $r,s,t \in \Z$, again with $h(x_0,x_2)$ unique up to
choice of $r$, $s$ and $t$.

(iii) Commutativity:

\begin{equation*}
\iota^{-1}_{1 \ 2} \langle v', (Id_V \otimes \co(x_2)) \co(x_1)v \rangle
=\iota^{-1}_{2 \ 1} \langle v', (T \otimes Id_V)
(Id_V \otimes \co(x_1)) \co(x_2)v \rangle.
\end{equation*}

(iv) Associativity:

\begin{equation*}
\iota^{-1}_{1 \ 2} \langle v', (Id_V \otimes \co(x_2)) \co(x_1)v \rangle
=\left. \left( \iota^{-1}_{2 \ 0}
\langle v', (\co(x_0) \otimes Id_V) \co(x_2)v \rangle
\right) \right|_{x_0=x_1-x_2}.
\end{equation*}
\end{proposition}

\begin{proposition} \label{P:RCAtoJI}
In the presence of the other VOC axioms, right and left
rationality, commutativity, and associativity imply the Jacobi
identity.
\end{proposition}

For the proofs of Propositions \ref{P:RR,LR,com,ass} and
\ref{P:RCAtoJI} see \cite{K} or \cite{K3}.

\section{An isomorphism between the categories of GVOCs and VOCs}
\label{S:iso1} We defined geometric vertex operator coalgebras to
motivate the definition of vertex operator coalgebras.  We are now
in a position to prove that these definitions actually define
isomorphic categories. The category of GVOCs of rank $d$ has the
set of all GVOCs of rank $d$ as its objects and its morphisms are
those linear maps between GVOCs that are $\mu$ invariant. The
category of VOCs of rank $d$ has all VOCs of rank $d$ as objects
and for morphisms, linear maps between VOCs which are coproduct
invariant as well as preserving the covacuum and co-Virasoro maps.

We will define a functor from GVOCs to VOCs and then a functor
from VOCs to GVOCs. Finally, we will show that these two functors
are inverses to each other, thus giving an isomorphism.

\subsection{A map from GVOCs to VOCs}
In this section we construct a map from the category of geometric
vertex operator coalgebras to the category of vertex operator
coalgebras.  In Section \ref{S:iso}, we will prove that, in fact,
this map is an isomorphism.

Let $(V,\mu)$ be a geometric vertex operator coalgebra. We define
a linear map $c^{\mu}: V \to \C$ by

\begin{equation}
c^{\mu}= \mu_0((1,\mathbf{0})),
\end{equation}

\noindent a linear map $\rho^{\mu}: V \to \C$ by

\begin{equation}
\rho^{\mu}= \left. -\frac{d}{d \ep} \mu_0(1,(0,\ep,0,0,\ldots))
\right|_{\ep=0},
\end{equation}

\noindent and a linear map $\co^{\mu} : V \mapsto (V \otimes
V)[[x,x^{-1}]]$ by

\begin{equation*}
\mathrm{Res}_x x^n \langle v',\co^{\mu}(x) v \rangle
=\mathrm{Res}_z z^n \langle v',\mu((z^{-1} ;\mathbf{0},
(1,\mathbf{0}), (1,\mathbf{0}) ))v \rangle,
\end{equation*}

\noindent for $v' \in (V \otimes V)', v \in V$ where
$\mathrm{Res}_x$ means taking the coefficient of the $x^{-1}$ term
in the given series and $\mathrm{Res}_z$ means taking the residue
of the function at the singularity $z=0$.
%(The residue exists by the meromorphicity axiom.)

\begin{proposition}
If the rank of $(V, \mu)$ is $d$, then the quadruple $(V,
\co^{\mu}, c^{\mu}, \rho^{\mu})$ is a vertex operator coalgebra of
rank $d$.
\end{proposition}

\begin{proof}
We use the GVOC axioms to prove the VOC axioms with respect to
$(V, \co^{\mu}, c^{\mu}, \rho^{\mu})$.

1. Left Counit:

\begin{align*}
\langle v',(c^{\mu} \otimes Id_V) \co^\mu(x)(v) \rangle |_{x=z} &
=
\left. \langle v',(\mu_0((1, \mathbf{0})) \otimes Id_V) \co^\mu(x)v \rangle \right|_{x=z} \\
& = \langle v', \mu_0((1, \mathbf{0})) \ _1*_{-1} \ \mu_2((z^{-1};
\mathbf{0}, (1,\mathbf{0}),(1,\mathbf{0}) ))v \rangle \\
& =  \langle v', \mu_1((1,\mathbf{0}) \ _1\infty_{-1} \
(z^{-1};\mathbf{0}, (1,\mathbf{0}),(1,\mathbf{0}) ))v \rangle \\
& = \langle v', \mu_1(\mathbf{0},(1,\mathbf{0}) )v \rangle \\
& = \langle v', v \rangle.
\end{align*}

(The second equality uses projection maps to show convergence and
the last step uses a trivial application of the grading axiom.)

2. Cocreation:

\begin{align*}
\langle v',(Id_V \otimes c^{\mu}) \co^\mu(x)(v) \rangle |_{x=z} &
= \left. \langle v',(Id_V \otimes
\mu_0((1,\mathbf{0}))) \co^\mu(x)v \rangle \right|_{x=z} \\
& =  \langle v', \mu_0((1,\mathbf{0})) \ _1*_{-2} \
\mu_2( (z^{-1};\mathbf{0}, (1,\mathbf{0}),(1,\mathbf{0}) ))v \rangle \\
& =  \langle v', \mu_1((1,\mathbf{0}) \ _1\infty_{-2} \
(z^{-1};\mathbf{0}, (1,\mathbf{0}),(1,\mathbf{0}) ))v \rangle \\
& =  \langle v', \mu_1(\mathbf{0}, (1,(-z,0,0,\ldots)) )v \rangle.
\end{align*}

By the meromorphicity axiom, $ \langle v', \mu_1(\cdot)v \rangle$
is a meromorphic function on $K^*(1)$, so $\langle v',
\mu_1(\mathbf{0}, (1,(-z,0,0,\ldots)) )v \rangle$ is a polynomial
in $z$.  Thus for all $v \in V$, $v' \in V'$, $\langle v',(Id_V
\otimes c^{\mu}) \co^{\mu}(x)(v) \rangle \in V[x] $. We can also
see that

\begin{eqnarray*}
\lim_{z \to 0} \langle v', \mu_1(\mathbf{0}, (1,(-z,0,0,\ldots))
)v \rangle & = &
    \langle v', \mu_1(\mathbf{0}, (1,\mathbf{0}) )v \rangle \\
        & = & \langle v', v \rangle.
\end{eqnarray*}

3. Truncation:

By the meromorphicity axiom, for any $v \in V$ there exists $N(v)
\in \Z_+$ such that for all $v' \in (V \otimes V)'$, the power of
$z$ in $\langle v',\mu_2((z^{-1};\mathbf{0},
(1,\mathbf{0}),(1,\mathbf{0}) ))v \rangle$ is less than $N(v)$.
But

\begin{equation*}
\langle v',\mu_2((z^{-1};\mathbf{0}, (1,\mathbf{0}),(1,\mathbf{0})
))v \rangle =\left. \langle v',\co^{\mu}(x)v \rangle
\right|_{x=z},
\end{equation*}

so the number of positive powers of x in $\co^{\mu}(x)v$ must be
less than $N(v)$ as well.

4. Jacobi Identity:

The idea behind proving the Jacobi identity is to prove right and
left rationality as well as commutativity and associativity.
Appealing to Proposition \ref{P:RCAtoJI}, this is equivalent to
the Jacobi identity.

We will start by obtaining right rationality and commutativity.

\begin{align*}
\langle v', (Id_V & \otimes \co^{\mu}(x_2))\co^{\mu}(x_1)v \rangle |_{x_i=z_{-i}} \\
&=\langle v', (Id_V \otimes \mu_2((z_{-2}^{-1};\mathbf{0},
(1,\mathbf{0}),(1,\mathbf{0}) ) ))
\mu_2((z_{-1}^{-1};\mathbf{0}, (1,\mathbf{0}),(1,\mathbf{0}) ) )v \rangle  \\
&=\langle v', \mu_3((z_{-2}^{-1};\mathbf{0},
(1,\mathbf{0}),(1,\mathbf{0})) \ _1\infty_{-2} \
(z_{-1}^{-1};\mathbf{0}, (1,\mathbf{0}),(1,\mathbf{0})) )v \rangle \\
&=\langle v', \mu_3((z^{-1}_{-1},z^{-1}_{-2}; \mathbf{0},
(1,\mathbf{0}), (1,\mathbf{0}), (1,\mathbf{0}) ))v \rangle
\end{align*}

\noindent for any $z_{-1},z_{-2} \in \C^{\times}$ for which this
sewing is well defined, i.e., for $|z_{-1}| > |z_{-2}|$.
Similarly,

\begin{multline*}
\langle v', (Id_V \otimes \co^{\mu}(x_1))\co^{\mu}(x_2)v \rangle
|_{x_i=z_{-i}} =\langle v', \mu_3((z^{-1}_{-2},z^{-1}_{-1};
\mathbf{0}, (1,\mathbf{0}), (1,\mathbf{0}), (1,\mathbf{0}) ))v
\rangle
\end{multline*}

\noindent for $|z_{-2}| > |z_{-1}|$.  By the meromorphicity axiom,

\begin{equation*}
\langle v', \mu_3((z^{-1}_{-1},z^{-1}_{-2}; \mathbf{0},
(1,\mathbf{0}), (1,\mathbf{0}), (1,\mathbf{0}) ))v \rangle
=\frac{g(z_{-1},z_{-2})}{z_{-1}^r z_{-2}^s (z_{-1} -z_{-2})^t}
\end{equation*}

\noindent where $g(z_{-1},z_{-2}) \in \C[z_{-1}, z_{-2}]$ and
$r,s,t \in \Z$.  Thus, since

\begin{equation*}
\langle v', (Id_V \otimes \co^{\mu}(x_2))\co^{\mu}(x_1)v \rangle
|_{x_i=z_{-i}} =\frac{g(z_{-1},z_{-2})}{z_{-1}^r z_{-2}^s (z_{-1}
-z_{-2})^t}
\end{equation*}

\noindent and $\langle v', (Id_V \otimes
\co^{\mu}(x_2))\co^{\mu}(x_1)v \rangle$ has only finitely many
positive powers of $x_1$ and hence only finitely many negative
powers of $x_2$ by evaluating weights on homogeneous vectors,

\begin{equation*}
\langle v', (Id_V \otimes \co^{\mu}(x_2))\co^{\mu}(x_1)v \rangle
=\iota_{1 \ 2} \frac{g(x_1,x_2)}{x_1^r x_2^s (x_1 -x_2)^t},
\end{equation*}

\noindent i.e., $V$ satisfies right rationality. By the
permutation axiom,

\begin{multline*}
\langle v', \mu_3((z^{-1}_{-1},z^{-1}_{-2};
(1,\mathbf{0}),\mathbf{0},  (1,\mathbf{0}), (1,\mathbf{0}) ))v \rangle \\
=\langle v', (T \otimes Id_V) \mu_3((z^{-1}_{-2},z^{-1}_{-1};
\mathbf{0}, (1,\mathbf{0}), (1,\mathbf{0}) ))v \rangle.
\end{multline*}

\noindent Thus

\begin{multline*}
\iota^{-1}_{1 \ 2} \langle v', (Id_V \otimes
\co^{\mu}(x_2))\co^{\mu}(x_1)v \rangle
=\frac{g(x_1,x_2)}{x_1^r x_2^s (x_1 -x_2)^t} \\
=\iota^{-1}_{2 \ 1} \langle v', (T \otimes Id_V) (Id_V \otimes
\co^{\mu}(x_1))\co^{\mu}(x_2)v \rangle
\end{multline*}

\noindent proving commutativity.  For left rationality and
associativity, we use the same technique to argue that

\begin{align*}
\langle v', &(\co^{\mu}(x_0) \otimes Id_V)\co^{\mu}(x_2)v \rangle
|_{x_0=(z_{-1}-z_{-2}),x_2=z_{-2}} \\
&=\langle v', \mu_2(((z_{-1}-z_{-2})^{-1}; \mathbf{0},
(1,\mathbf{0}), (1,\mathbf{0}) ))
\ _1*_{-2} \ \mu_2((z_{-2}^{-1}; \mathbf{0}, (1,\mathbf{0}), (1,\mathbf{0}) ))v \rangle\\
&=\langle v', \mu_3((z_{-1},z_{-2}; \mathbf{0}, (1,\mathbf{0}),
(1,\mathbf{0}), (1,\mathbf{0}) ))v \rangle.
\end{align*}

\noindent for $|z_{-2}| > |z_{-1}-z_{-2}|$. Recall that the
right-hand side of this equation is equal to

\begin{equation*}
\frac{g(z_{-1},z_{-2})}{z_{-1}^r z_{-2}^s (z_{-1} -z_{-2})^t}
=\frac{h(z_{-1}-z_{-2},z_{-2})} {((z_{-1}-z_{-2})+z_{-2})^r
z_{-2}^s (z_{-1} -z_{-2})^t}
\end{equation*}

\noindent for some $h(z_{-1}-z_{-2},z_{-2}) \in \C[z_{-1},
z_{-2}]$.  Thus as above,

\begin{equation*}
\langle v', (Id_V (\co^{\mu}(x_0) \otimes Id_V)\co^{\mu}(x_2)v
\rangle =\iota_{2 \ 0} \frac{h(x_0,x_2)}{(x_0+x_2)^r x_2^s
x_0^{t}},
\end{equation*}

\noindent that is, left rationality holds.  In addition, we see
that

\begin{align*}
 \left( \iota^{-1}_{2 \ 0} \langle v', (Id_V (\co^{\mu}(x_0)
\otimes Id_V)\co^{\mu} ( x_2)v \rangle
\right) & |_{x_0=x_1-x_2} \\
&=\frac{h(x_{-1}-x_{-2},x_{-2})}
{((x_{-1}-x_{-2})+x_{-2})^r x_{-2}^s (x_{-1} -x_{-2})^t} \\
&=\iota^{-1}_{1 \ 2} \langle v', (Id_V \otimes
\co^{\mu}(x_2))\co^{\mu}(x_1)v \rangle,
\end{align*}

\noindent which is, of course, associativity.

5. Virasoro algebra:

In order to examine the Virasoro algebra structure of the $L(k)$
operators, we must first determine what

\begin{equation*}
(\rho \otimes Id_V) \co(x) = \sum_{k \in \Z} L(k) x^{k-2}
\end{equation*}

\noindent forces the $L(k)$ operators to be. We will make use of
the linear functionals $\mathcal{L}_I(z)$ described in Section
\ref{S:L_I}. If we choose $v' \in V'$ and $v \in V$ then apply
$\mathcal{L}_I(z)$ to $\langle v', \mu_1(\cdot)v \rangle$, there
are two different ways to interpret the result.  On the one hand,

\begin{align*}
\mathcal{L}_I(z) &\langle v', \mu_1(\cdot)v \rangle \\
&=\left.\left( \frac{d}{d \ep} \langle v',
\mu_1(((1,(0,-\ep,0,0,\ldots))) \ _1\infty_{-1} \
(z^{-1};\mathbf{0}, (1,\mathbf{0}),(1,\mathbf{0}) )) v \rangle
\right) \right|_{\ep=0} \\
&= \left. \left( \frac{d}{d \ep} \langle v', \mu_0(\mathbf{0},
(1,(0,-\ep,0,0,\ldots)) ) \ _1*_{-1} \ \co^\mu(x)) v \rangle
\right)  \right|_{\ep=0, x=z} \\
&=  \left. \langle v',\frac{-d}{d \ep}
\mu_0(\mathbf{0},(1,(0,\ep,0,0,\ldots)) )v \rangle \right|_{\ep=0}
\left._1*_{-1}  \langle v',\co^\mu(x)) v \rangle \right|_{x=z} \\
&= \left. \langle v', (\rho^{\mu} \otimes Id_V) \co^{\mu}(x) v
\rangle \right|_{x=z}
\end{align*}

On the other hand, viewing $\mathcal{L}_I(z)$ as an element of the
tangent space of $K^*(1)$ at $(\mathbf{0}, (1,\mathbf{0}))$, Lemma
\ref{L:L_I} implies

\begin{align*}
\mathcal{L}_I(z)&\langle v', \mu_1(\cdot)v \rangle \\
=&  \sum_{k \in \Z_+} \left. z^{-k-2} \frac{-\partial}{\partial
A_k^{(-1)}} \langle v', \mu_1((a_0^{(1)},A^{(1)}),A^{(-1)})v
\rangle
\right|_{A^{(-1)},A^{(1)}= \mathbf{0}, \ a_0^{(1)}=1}  \\
&+\left. z^{-2} \frac{-\partial}{\partial a_0^{(1)}} \langle v',
\mu_1((a_0^{(1)},A^{(1)}),A^{(-1)})v \rangle
\right|_{A^{(-1)},A^{(1)}= \mathbf{0}, \ a_0^{(1)}=1} \\
&+ \sum_{k \in \Z_+} \left. z^{k-2} \frac{-\partial}{\partial
A_k^{(1)}} \langle v', \mu_1((a_0^{(1)},A^{(1)}),A^{(-1)})v
\rangle \right|_{A^{(-1)},A^{(1)}= \mathbf{0}, \ a_0^{(1)}=1} \\
=& \sum_{k \in \Z_+} \langle v',z^{-k-2} L(-k)v \rangle + \langle
v',z^{-2} L(0) v \rangle +\sum_{k \in \Z_+} \langle v', z^{k-2}
L(k)v \rangle \\
=& \left. \langle v', \sum_{k \in \Z}  L(k) x ^{k-2}v \rangle
\right|_{x=z}
\end{align*}

\noindent where the $L(k)$'s are defined by

\begin{equation*}
\langle v',L(k) v \rangle =\left.-\frac{\partial}{\partial
A_k^{(1)}} \langle v', \mu_1((a_0^{(1)},A^{(1)}),A^{(-1)})v
\rangle \right|_{A^{(-1)},A^{(1)}= \mathbf{0}, \ a_0^{(1)}=1},
\end{equation*}
\begin{equation*}
\langle v',L(0) v \rangle =\left.-\frac{\partial}{\partial
a_0^{(1)}}\langle v', \mu_1((a_0^{(1)},A^{(1)}),A^{(-1)})v \rangle
\right|_{A^{(-1)},A^{(1)}= \mathbf{0}, \ a_0^{(1)}=1},
\end{equation*}
\begin{equation*}
\langle v',L(-k) v \rangle =\left.-\frac{\partial}{\partial
A_k^{(-1)}} \langle v', \mu_1((a_0^{(1)},A^{(1)}),A^{(-1)})v
\rangle \right|_{A^{(-1)},A^{(1)}= \mathbf{0}, \ a_0^{(1)}=1}
\end{equation*}

\noindent for $k>0$.

This geometric definition of the $L(k)$ operators is identical to
the definition of the $L(k)$ operators in the vertex operator
algebra setting of \cite{H}. Thus we may claim the proof from
vertex operator algebras (\cite{H}, (5.4.26)) that these operators
satisfy the Virasoro relations with the bracket as defined in
\cite{H}.

6. Grading:  If we assume that $v \in V_{(k)}$, then observe that
by the definition of $L(0)$ and the GVOC grading axiom

\begin{eqnarray*}
\langle v',L(0)v \rangle &= & \left.\left(
-\frac{\partial}{\partial a_0^{(1)}}\langle v', \mu_1(\mathbf{0},
(a_0,\mathbf{0}))v \rangle \right)
\right|_{ a_0=1} \\
&= & \left.\left( -\frac{\partial}{\partial a_0^{(1)}}
(a_0^{(1)})^{-k} \langle v',v \rangle
\right) \right|_{ a_0=1} \\
&= & \langle v',k v \rangle.
\end{eqnarray*}

7. $L(1)$-Derivative:

\begin{align*}
&\left. \langle v',(L(1) \otimes Id_V) \co^\mu(x)v \rangle \right|_{x=z} \\
&= \left. \left( \frac{-\partial}{\partial z_0}\langle v',
(\mu_1(\mathbf{0}, (1,(z_0,0,0,\ldots)) ) \otimes Id_V)\co^\mu(x)v
\rangle \right)
\right|_{x=z, z_0=0} \\
&= \left. \left( \frac{-\partial}{\partial z_0}\langle v',
\mu_1(\mathbf{0}, (1,(z_0,0,0,\ldots)) )
_1*_{-1}\mu_2((z^{-1};\mathbf{0}, (1,\mathbf{0}),
(1,\mathbf{0}) ))v \rangle \right) \right|_{z_0=0} \\
&= \left. \left( \frac{-\partial}{\partial z_0}\langle v',
\mu_2(\mathbf{0}, (1,(z_0,0,0,\ldots)) )
_1\infty_{-1}(z^{-1};\mathbf{0}, (1,\mathbf{0}),
(1,\mathbf{0}) ))v \rangle \right) \right|_{z_0=0} \\
&= \left. \left( \frac{-\partial}{\partial z_0} \langle v',
\mu_2(z-z_0; \mathbf{0}, (1,\mathbf{0}), (1,\mathbf{0}) )v \rangle
\right)
\right|_{z_0=0} \\
&= \left. \left( \frac{-\partial}{\partial z_0} \left( \left.
\langle v', \co^{\mu}(x)v \rangle \right|_{x=z-z_0} \right)
\right) \right|_{z_0=0} \\
&=  \left. \left( \left. \langle v', \frac{-\partial}{\partial x}
\co^{\mu}(x)v \rangle \right|_{x= z-z_0} \frac{\partial}{\partial
z_0}(z-z_0)
\right) \right|_{z_0=0} \\
&=  \left. \langle v', \frac{\partial}{\partial x} \co^{\mu}(x)v
\rangle \right|_{x= z}
\end{align*}

\noindent This completes the proof that all axioms hold for $(V,
\co^{\mu}, c^{\mu}, \rho^{\mu})$.

\end{proof}

\subsection{A map from VOCs to GVOCs} \label{S:VOC_GVOC}
We now construct a map from the category of vertex operator
coalgebras to the category of geometric vertex operator
coalgebras.

Let $(V, \co, c, \rho)$ be a VOC of rank $d$. Let $\mu^{\co}_n:
K^*(n) \mapsto \mathcal{H}_V(1,n)$ be defined by

\begin{multline} \label{E:def2}
\langle v', \mu_n^{\co}(z^{-1}_{-n+1},\ldots,z^{-1}_{-1};A^{(-n)},
(a_0^{(-n+1)},A^{(-n+1)}),
\ldots,\\
(a_0^{(-1)},A^{(-1)}),(a_0^{(1)},A^{(1)}) )v \rangle \\
=\iota^{-1}_{1 \cdots n-1} \langle v', \left( (a_0^{(-1)})^{-L(0)}
e^{-\sum_{j \in \Z_+} A^{(-1)}_j L(-j)} \otimes \cdots \otimes
(a_0^{(-n+1)})^{-L(0)}
\right. \\
\left.  e^{-\sum_{j \in \Z_+} A^{(-n+1)}_j L(-j)} \otimes
e^{-\sum_{j \in \Z_+} A^{(-n)}_j L(-j)} \right)
(\underbrace{Id_V \otimes \cdots \otimes Id_V}_{n-2} \otimes \co(x_{n-1})) \cdots \\
(Id_V \otimes \co(x_{2})) \co(x_{1}) \left. e^{-\sum_{j \in \Z_+}
A^{(1)}_j L(j)} ( a_0^{(1)})^{-L(0)}v \rangle
\right|_{x_{i}=z_{-i}}
\end{multline}

\noindent for $n \in \Z_+$ and

\begin{equation*}
\mu_0^{\co}((1, A^{(1)}))v =  c e^{-\sum_{j \in \Z_+} A^{(1)}_j
L(j)} v
\end{equation*}

\noindent for all $v \in V$, $v' \in (V^{\otimes n})'$.  The
right-hand side of (\ref{E:def2}) is well defined since $v$ and
$v'$ are non-zero in on finitely many weights and degree raising
operators are to the left of degree lowering operators.

At this point we need a lemma that gives us some information about
absolute convergence when two sewings are composed.  In and of
itself this lemma says little because of its strong conditions,
but it will be necessary in the general proof that $(V, \{
\mu^{\co}_n \}_{n \in \Z})$ is a GVOC.

\begin{lemma} \label{L:doubleconvergence}
Let $\ell,m \in \N$ and $n \in \Z_+$ such that $m+n > 1$.  Choose
$Q_1 \in K^*(\ell), \ Q_2 \in K^*(m), \ Q_3 \in K^*(n)$ and
integers $1 \leq i \leq m$ and $1 \leq j \leq n$.  Assume that the
sewing $(Q_1 \ _1\infty_{-i} \ Q_2) \ _1\infty_{-j} Q_3$ exists.
Then if, for all $v' \in (V^{\otimes \ell+m+n-2})'$, $v \in V$,

\begin{equation*}
\langle v', \mu_{\ell+m-1}^{\co}(Q_1 \ _1\infty_{-i} \ Q_2) v
\rangle = \langle v', (\mu_{\ell}^{\co}(Q_1) \ _1*_{-i} \
\mu_{m}^{\co}(Q_2)) v \rangle e^{\Gamma_1 d}
\end{equation*}

\noindent and

\begin{multline*}
\langle v', \mu_{\ell+m+n-2}^{\co}((Q_1 \ _1\infty_{-i} \ Q_2) \ _1\infty_{-j} \ Q_3) v \rangle = \\
\langle v',(\mu_{\ell+m-1}^{\co}(Q_1 \ _1 \infty_{-i} \ Q_2)  \
_1*_{-j} \ \mu_n^{\co}(Q_3) ) v \rangle e^{\Gamma_2 d}
\end{multline*}

\noindent for $\Gamma_1, \Gamma_2$ appropriately chosen, we have

\begin{multline} \label{E:226}
\langle v', \mu_{\ell+m+n-2}^{\co}((Q_1 \ _1\infty_{-i} \ Q_2) \ _1\infty_{-j} \ Q_3) v \rangle = \\
\langle v',((\mu_{\ell}^{\co}(Q_1) \ _1*_{-i} \ \mu_m^{\co}(Q_2))
\ _1*_{-j} \ \mu_n^{\co}(Q_3) ) v \rangle e^{(\Gamma_1+\Gamma_2)d}
\end{multline}

\noindent and, in particular, the right-hand side of Equation
\ref{E:226} does exist.

Similarly, for $1 \leq i \leq m+n-1$ and $1 \leq j \leq n$, assume
that the sewing $Q_1 \ _1\infty_{-i} \ (Q_2 \ _1\infty_{-j} Q_3)$
exists.  If, for all $v' \in (V^{\otimes \ell+m+n-2})'$, $v \in
V$,

\begin{equation*}
\langle v', \mu_{m+n-1}^{\co}(Q_2 \ _1\infty_{-j} \ Q_3) v \rangle \\
=\langle v', (\mu_{m}^{\co}(Q_2) \ _1*_{-j} \ \mu_{n}^{\co}(Q_3))
v \rangle e^{\Gamma_1 d}
\end{equation*}

\noindent and

\begin{multline*}
\langle v', \mu_{\ell+m+n-2}^{\co}(Q_1 \ _1\infty_{-i} \ (Q_2 \ _1\infty_{-j} Q_3)) v \rangle \\
=\langle v', (\mu_{\ell}^{\co}(Q_1) \ _1*_{-i} \
\mu_{m+n-1}^{\co}(Q_2  \ _1\infty_{-j} \ Q_3) ) v \rangle
e^{\Gamma_2 d}
\end{multline*}

\noindent for $\Gamma_1, \Gamma_2$ appropriately chosen, we have

\begin{multline} \label{E:226b}
\langle v', \mu_{\ell+m+n-2}^{\co}(Q_1 \ _1\infty_{-i} \ (Q_2 \ _1\infty_{-j} Q_3)) v \rangle \\
=\langle v', (\mu_{\ell}^{\co}(Q_1) \ _1*_{-i} \ (\mu_m^{\co}(Q_2)
\ _1*_{-j} \ \mu_n^{\co}(Q_3) )) v \rangle
e^{(\Gamma_1+\Gamma_2)d}
\end{multline}

\noindent and, in particular, the right-hand side of Equation
\ref{E:226b} does exist.
\end{lemma}

\begin{proof}
This result uses the Fischer-Grauert Theorem (Theorem 3.4.3 in
\cite{H}) and is essentially the double absolute convergence
result given in part (5e) of Proposition 5.4.1 in \cite{H}.  (It
also follows Proposition 7.1 in \cite{B}).  To prove the first
half of the lemma, let $t_1, t_2 \in \C^{\times}$, then multiply
the local coordinate map at the $-i$-th puncture of $Q_2$ by $t_1$
and at the $-j$-th puncture of $Q_3$ by $t_2$ to obtain $Q_2(t_1)
\in K^*(m)$ and $Q_3(t_2) \in K^*(n)$, respectively.  There exists
a neighborhood of $(1,1) \in \C^{\times} \times \C^{\times}$ such
that the sewing still exists.  Using the definitions necessary to
expand both sides of (\ref{E:226}), we can observe that
substituting $Q_2(t_1)$ and $Q_3(t_2)$ for $Q_2$ and $Q_3$ on the
left-hand side of (\ref{E:226}) shows that both sides are equal
when viewed as formal power series in $t_1$ and $t_2$.  (This
substitution is valid because of the Fischer-Grauert Theorem.)
Since the left-hand side is doubly absolutely convergent at
$t_1=t_2=1$, so is the right-hand side.

The argument for the second half of the lemma is similar.
\end{proof}

\noindent We now move to the main proposition of this section.

\begin{proposition} \label{P:VOCtoGVOC}
The pair $(V, \{ \mu^{\co}_n \}_{n \in \Z})$ is a geometric
 vertex operator coalgebra of rank $d$, where $d$ is the rank of the VOC
$(V, \co, c, \rho)$.
\end{proposition}

\begin{proof}
Essentially, we use the axioms of a VOC to obtain the axioms of a
GVOC:

1. Grading: Let $v' \in V'$ and $v \in V_{(k)}$. Then by the VOC
grading axiom (\ref{E:VOCgrading}),

\begin{eqnarray*}
\langle v', \mu_1^{\co}(\mathbf{0}, (a,\mathbf{0}) )v \rangle
&=&\langle v', a^{-L(0)}v \rangle \\
&=&a^{-k} \langle v', v \rangle.
\end{eqnarray*}

2. Meromorphicity:  From right rationality (Proposition
\ref{P:RR,LR,com,ass}(i)) and the definition of $\mu_n^{\co}$, we
see that the map $Q \to \langle v', \mu_n^{\co}(Q)v \rangle$ is
meromorphic on $K(n)$. Fixing $v \in V$, we must show that for
some $N(v) \in \Z_+$ the degree of $z_{-i}$ is less than $N(v)$ in
$\langle v', \mu_n^{\co}(Q)v \rangle$ for any $1 \leq i \leq n-1$
and $v' \in (V^{\otimes n})'$. We may assume $Q=(z^{-1}_{-n+1},
\ldots, z^{-1}_{-1}; \mathbf{0}, (1,\mathbf{0}), \ldots,
(1,\mathbf{0}) )$.

Considering repeated application of commutativity (Propositions
\ref{P:RR,LR,com,ass} (iii)), we see that

\begin{multline*}
\langle v', \mu_n^{\co}(Q)v \rangle \\
=\iota^{-1}_{1 \cdots n-1} \langle v', (\underbrace{Id_V \otimes
\cdots \otimes Id_V}_{n-2} \otimes \co(x_{n-1})) \cdots (Id_V
\otimes \co(x_{2})) \co(x_{1}) \left.  v \rangle
\right|_{x_{\ell}=z_{-\ell}}
\end{multline*}
\begin{multline*}
=\iota^{-1}_{1 \cdots i-2 \ i \ i-1 \ i+1 \cdots n-1} \langle v',
(\underbrace{Id_V \otimes \cdots \otimes Id_V}_{n-2} \otimes \co(x_{n-1})) \cdots \\
(\underbrace{Id_V \otimes \cdots \otimes Id_V}_{i} \otimes
\co(x_{i+1})) (\underbrace{Id_V \otimes \cdots \otimes Id_V}_{i-2}
\otimes
(T \otimes Id_V)(Id_V \otimes \co(x_{i-1}))\co(x_i)) \\
(\underbrace{Id_V \otimes \cdots \otimes Id_V}_{i-3} \otimes
\co(x_{i-2})) \cdots (Id_V \otimes \co(x_{2})) \co(x_{1}) \left. v
\rangle \right|_{x_{\ell}=z_{-\ell}}
\end{multline*}
\begin{multline*}
=\iota^{-1}_{i \ 1 \cdots \widehat{i} \cdots n-1} \langle v',
(\underbrace{Id_V \otimes \cdots \otimes Id_V}_{n-2} \otimes \co(x_{n-1})) \cdots \\
(\underbrace{Id_V \otimes \cdots \otimes Id_V}_{i} \otimes
\co(x_{i+1}))
(\underbrace{Id_V \otimes \cdots \otimes Id_V}_{i-2} \otimes \\
(T \otimes Id_V)(Id_V \otimes \co(x_{i-1})) \cdots (T \otimes
Id_V)(Id_V \otimes \co(x_{1})) \co(x_{i}) \left.  v \rangle
\right|_{x_{\ell}=z_{-\ell}}.
\end{multline*}

\noindent By the truncation axiom, there exists $N(v) \in \Z_+$
such that the term $\co(x_i)v$ has less than $N(v)$ positive
powers of $x_i$.  Thus we have less that $N(v)$ positive powers of
$z_{-i}$.

3. Permutation: Given $Q \in K^*(n)$, we note that $S_n$ is
generated by the actions of the transpositions $\sigma=(i \ i+1)$
 for $i = 1,2, \ldots, n-1$ so it suffices to show that the axiom holds for these.
If $i<n-1$, the definition of the action and Proposition
\ref{P:RR,LR,com,ass} (commutativity) suffice.  The calculation is
somewhat lengthy and found in \cite{K}. If $i=n-1$, the argument
requires not commutativity but the algebraic facts
(\ref{E:shift}), (\ref{E:distrib}), (\ref{E:expL(0)}) and
(\ref{E:antisym}).  Again, see \cite{K} for the explicit
calculation.

4. Sewing:

To prove the sewing axiom, consider $Q_1 \in K^*(m)$ and $Q_2 \in
K^*(n)$ and $1 \leq i \leq n$.  Our proof will use the approach in
\cite{B} and consider 8 steps. Step (a) will be the case that
$m=0,1$, $n=1$, (b) - (f) will cover special cases needed for the
induction, (g) is induction on $m$, and (h) is induction on $n$.

Throughout these steps, we will make use of the basis $\{
e_{l^{(k)}}^{(k)} |l^{(k)}=1, \ldots,\dim V_{(k)}; \\ k \in \Z \}$
of $V^{\otimes m}$ and $\{ (e_{l^{(k)}}^{(k)})^* |l^{(k)}=1,
\ldots,\dim V_{(k)}; \ k \in \Z \}$ the corresponding dual basis.

Step (a): Since $K(1)=K^*(1)$, the case $Q_1, Q_2 \in K^*(1)=K(1)$
where \\
$Q_1 \ _1\infty_{-1} \ Q_2$ exists is proven in \cite{H}.
Thus we have

\begin{equation*}
\langle v', \mu_{1}^{\co} (Q_1 \ _1\infty_{-1} \ Q_2)v \rangle =
\langle v',(\mu_1^{\co}(Q_1) _1*_{-1} \mu_1^{\co}(Q_2) v \rangle
e^{-\Gamma(A^{(1)}, B^{(-1)},a_0^{(1)})d}
\end{equation*}

\noindent The case $Q_1 \in K^*(0)$, $Q_2 \in K^*(1)$ is similar,
but without needing $v' \in V'$.

Step (b): Now consider the case $i= 2$ with

\begin{equation*}
Q_1=(A^{(-1)}, (a_0^{(1)},A^{(1)}) ) \in K^*(1),
\end{equation*}
\begin{equation*}
Q_2=(\zeta^{-1};\mathbf{0}, (b_0^{(-1)},B^{(-1)}), (1,\mathbf{0})
) \in K^*(2),
\end{equation*}

\noindent such that $Q_1 \ _1\infty_{-i} \ Q_2$ exists.

\begin{equation*}
\langle v', (\mu_1^{\co}(Q_1) _1*_{-i} \mu_2^{\co}(Q_2))_t v\rangle \\
\end{equation*}
\begin{multline*}
=\sum_{k \in \Z} \sum_{l^{(k)}=1}^{\dim V_{(k)}} \langle v',
\text{\Huge ( }Id_V \otimes \left( e^{-\sum_{j \in \Z_+}
A^{(-1)}_j L(-j)}
 e^{-\sum_{j \in \Z_+} A^{(1)}_j L(j)} ( a_0^{(1)})^{-L(0)} \right)  \\
 e^{(k)}_{l^{(k)}} (e^{(k)}_{l^{(k)}})^* \text{\Huge )}
  \left(
(b_0^{(-1)})^{-L(0)} e^{-\sum_{j \in \Z_+} B^{(-1)}_j L(-j)}
\otimes Id_V \right) \left. \co(x) v \rangle \right|_{x=\zeta} t^k
\end{multline*}
\begin{multline*}
=\sum_{k \in \Z} \sum_{l^{(k)}=1}^{\dim V_{(k)}} \langle v',
\text{\Huge (} Id_V \otimes \left( e^{-\sum_{j \in \Z_+}
A^{(-1)}_j L(-j)}
 e^{-\sum_{j \in \Z_+} A^{(1)}_j L(j)} ( a_0^{(1)})^{-L(0)} t^{L(0)} \right)  \\
 e^{(k)}_{l^{(k)}} (e^{(k)}_{l^{(k)}})^* \text{\Huge )}
  \left(
(b_0^{(-1)})^{-L(0)} e^{-\sum_{j \in \Z_+} B^{(-1)}_j L(-j)}
\otimes Id_V \right) \co(x) \left. v \rangle \right|_{x=\zeta}
\end{multline*}
\begin{multline*}
=\langle v', \left(Id_V \otimes e^{-\sum_{j \in \Z_+} A^{(-1)}_j
L(-j)}
 e^{-\sum_{j \in \Z_+} A^{(1)}_j L(j)} ( a_0^{(1)} t^{-1})^{-L(0)} \right) \\
  \left( (b_0^{(-1)})^{-L(0)} e^{-\sum_{j \in \Z_+} B^{(-1)}_j L(-j)}  \otimes Id_V \right)
\co(x) \left. v \rangle \right|_{x=\zeta}
\end{multline*}
\begin{multline*}
=\langle v', \left((b_0^{(-1)})^{-L(0)} e^{-\sum_{j \in \Z_+}
B^{(-1)}_j L(-j)}
\otimes e^{-\sum_{j \in \Z_+} A^{(-1)}_j L(-j)} \right )\\
\left(Id_V \otimes  e^{-\sum_{j \in \Z_+} A^{(1)}_j L(j)} (
a_0^{(1)} t^{-1})^{-L(0)}\right)
\co(x) ( a_0^{(1)} t^{-1})^{L(0)}  e^{\sum_{j \in \Z_+} A^{(1)}_j L(j)} \\
e^{-\sum_{j \in \Z_+} A^{(1)}_j L(j)} ( a_0^{(1)} t^{-1})^{-L(0)}
\left. v \rangle \right|_{x=\zeta}.
\end{multline*}

\noindent On the other hand, if we consider $Q_2$ as above and

\begin{equation*}
Q_1(t)=(A^{(-1)}, (a_0^{(1)}t^{-1},A^{(1)}) ) \in K^*(1),
\end{equation*}

\noindent then by Example \ref{E:sewing5}, we observe that the
image under $\mu_2^{\co}$ of the sewing $Q_1(t) \ _1\infty_{-2} \
Q_2$ exists for $t=1$.  Thus

\begin{multline*}
\mu_2^{\co}(Q_1(t) \ _1\infty_{-2} \ Q_2)\\
=\langle v', \left((b_0^{(-1)})^{-L(0)} e^{-\sum_{j \in \Z_+}
B^{(-1)}_j L(-j)}
(a_0^{(1)}t^{-1})^{L(0)}e^{- \Theta_0^{(2)}L(0)} \right. \\
\left. e^{-\sum_{j \in \Z_+} (a_0^{(1)}t^{-1})^{-j}
\Theta_j^{(2)}L(-j)}
\otimes e^{-\sum_{j \in \Z_+} A^{(-1)}_j L(-j)} \right )\\
\co(a_0^{(1)} t^{-1} \hat{f}_2(x)) e^{-\sum_{j \in \Z_+} A^{(1)}_j
L(j)} ( a_0^{(1)} t^{-1})^{-L(0)} \left. v \rangle
\right|_{x=\zeta},
\end{multline*}

\noindent where $\Theta_j^{(2)}$, for $j \in \N$, and
$\hat{f}_2(x)$ are considered as functions of $a_0^{(1)}$,
$A^{(1)}$ and $x$.  It therefore remains to show that

\begin{multline} \label{E:1029}
\left(Id_V \otimes  e^{-\sum_{j \in \Z_+} A^{(1)}_j L(j)} (
a_0^{(1)} t^{-1})^{-L(0)}\right)
\co(x) ( a_0^{(1)} t^{-1})^{L(0)}  e^{\sum_{j \in \Z_+} A^{(1)}_j L(j)} \\
= \left( (a_0^{(1)} t^{-1})^{L(0)} e^{- \Theta_0^{(2)} L(0)}
e^{-\sum_{j \in \Z_+} (a_0^{(1)} t^{-1})^{-j} \Theta_j^{(2)}L(-j)}
\otimes Id_V \right ) \co( a_0^{(1)} t^{-1} \hat{f}_2(x)).
\end{multline}

\noindent Using (\ref{E:shift}), (\ref{E:L(0)L(k)}), and
Proposition \ref{P:sew_id3b}, we observe that the right-hand side
of (\ref{E:1029}) is equal to

\begin{multline*}
\left( (a_0^{(1)} t^{-1})^{L(0)}e^{- \Theta_0^{(2)}L(0)}
e^{-\sum_{j \in \Z_+} (a_0^{(1)}t^{-1})^{-j} \Theta_j^{(2)}L(-j)}
e^{a_0^{(1)} t^{-1}(-x+\hat{f}_2(x))L(1)}
\otimes Id_V \right) \\
\co(a_0^{(1)} t^{-1} x)
\end{multline*}
\begin{multline*}
=\left( e^{- \Theta_0^{(2)}L(0)} e^{-\sum_{j \in \Z_+}
\Theta_j^{(2)}L(-j)} e^{(-x+\hat{f}_2(x))L(1)} (a_0^{(1)}
t^{-1})^{L(0)} \otimes Id_V \right) \co(a_0^{(1)} t^{-1} x)
\end{multline*}
\begin{multline*}
=\left(e^{\sum_{k=-1}^{\infty} \left( \sum_{j \in \Z_+}
(a_0^{(1)}t^{-1})^{-j} A^{(1)}_j { -j+1 \choose k+1}x^{-j-k}
\right)L(-k)} (a_0^{(1)} t^{-1})^{L(0)} \otimes Id_V \right)
\co(a_0^{(1)} t^{-1} x)
\end{multline*}
\begin{multline*}
=\left( (a_0^{(1)} t^{-1})^{L(0)} \otimes Id_V \right)
\left(e^{\sum_{k=-1}^{\infty} \left( \sum_{j \in \Z_+}
(a_0^{(1)}t^{-1})^{-j-k} A^{(1)}_j
{ -j+1 \choose k+1}x^{-j-k} \right)L(-k)} \otimes Id_V \right) \\
\co(a_0^{(1)} t^{-1} x).
\end{multline*}

\noindent On the other hand, using (\ref{E:expL(0)}), we see that
the left-hand side of (\ref{E:1029}) is equal to

\begin{equation*}
\left(Id_V \otimes  e^{-\sum_{j \in \Z_+} A^{(1)}_j L(j)}\right)
\left(( a_0^{(1)} t^{-1})^{L(0)} \otimes Id_V \right) \co(
a_0^{(1)} t^{-1} x) e^{\sum_{j \in \Z_+} A^{(1)}_j L(j)}
\end{equation*}
\begin{equation*}
=\left(( a_0^{(1)} t^{-1})^{L(0)} \otimes Id_V \right) \left(Id_V
\otimes  e^{-\sum_{j \in \Z_+} A^{(1)}_j L(j)}\right) \co(
a_0^{(1)} t^{-1} x) e^{\sum_{j \in \Z_+} A^{(1)}_j L(j)}.
\end{equation*}

\noindent Thus it suffices to show that

\begin{multline} \label{E:1030}
\left(\sum_{k=-1}^{\infty} \left( \sum_{j \in \Z_+}
(a_0^{(1)}t^{-1})^{-j-k} A^{(1)}_j { -j+1 \choose k+1}x^{-j-k}
\right)L(-k) \otimes Id_V \right)
\co(a_0^{(1)} t^{-1} x) \\
=\co( a_0^{(1)} t^{-1} x) \sum_{j \in \Z_+} A^{(1)}_j L(j) -
\left(Id_V \otimes  \sum_{j \in \Z_+} A^{(1)}_j L(j)\right) \co(
a_0^{(1)} t^{-1} x).
\end{multline}

\noindent Define

\begin{equation*}
h(x_1) = \sum_{j \in \Z_+} A_j^{(1)} x_1^{-j+1}
\end{equation*}

\noindent and notice that

\begin{equation*}
\mathrm{Res}_{x_1} \left( h(x_1) \sum_{m \in \Z} L(m) x_1^{m-2}
\right) = \sum_{j \in \Z_+} A^{(1)}_j L(j).
\end{equation*}

Now recalling (\ref{E:rho_on_jac}) and observing that $h(x_1)$
commutes with all other terms, we observe that the right-hand side
of (\ref{E:1030}) is

\begin{multline*}
\mathrm{Res}_{x_1} \left( \co( a_0^{(1)} t^{-1} x)
h(x_1) \sum_{m \in \Z} L(m) x_1^{m-2} \right.\\
- \left(Id_V \otimes h(x_1) \sum_{m \in \Z} L(m) x_1^{m-2}\right)
\left. \co( a_0^{(1)} t^{-1} x) \right)
\end{multline*}
\begin{multline*}
=\mathrm{Res}_{x_1} \mathrm{Res}_{x_0} (a_0^{(1)}t^{-1}x)^{-1}
\delta \left( \frac{x_1-x_0}{a_0^{(1)}t^{-1}x} \right)
\left( h(x_1) \sum_{m \in \Z} L(m) x_0^{m-2} \otimes Id_V \right) \\
\co( a_0^{(1)} t^{-1} x)
\end{multline*}
\begin{multline*}
=\mathrm{Res}_{x_1} \mathrm{Res}_{x_0} x_1^{-1} \delta \left(
\frac{a_0^{(1)}t^{-1}x+x_0}{x_1} \right) \left( \sum_{j \in \Z_+}
A_j^{(1)} x_1^{-j+1}
 \sum_{m \in \Z} L(m) x_0^{m-2} \otimes Id_V \right) \\
\co( a_0^{(1)} t^{-1} x)
\end{multline*}
\begin{multline*}
=\mathrm{Res}_{x_0} \left( \sum_{j \in \Z_+} A_j^{(1)}
(a_0^{(1)}t^{-1}x +x)^{-j+1}
 \sum_{m \in \Z} L(m) x_0^{m-2} \otimes Id_V \right)
\co( a_0^{(1)} t^{-1} x) \\
\end{multline*}
\begin{multline*}
=\mathrm{Res}_{x_0} \left( \sum_{j \in \Z_+} A_j^{(1)}
\sum_{k=-1}^{\infty} {-j+1 \choose k+1} (a_0^{(1)}t^{-1}x)^{-j-k}
x_0^{k+1}
 \sum_{m \in \Z} L(m) x_0^{m-2} \otimes Id_V \right) \\
\co( a_0^{(1)} t^{-1} x).
\end{multline*}

\noindent But this is equal to the left-hand side of
(\ref{E:1030}). Therefore, since $\mu_2^{\co}(Q_1(t) \
_1\infty_{-2} \ Q_2)$ exists when $t=1$,

\begin{align*}
\mu_2^{\co}(Q_1(1) \ _1\infty_{-2} \ Q_2)
&=\langle v', (\mu_1^{\co}(Q_1) \ _1*_{-2} \ \mu_2^{\co}(Q_2))_t \left. v\rangle \right|_{t=1} \\
&=\langle v', (\mu_1^{\co}(Q_1) \ _1*_{-2} \  \mu_2^{\co}(Q_2))
v\rangle.
\end{align*}

Step (c): Let $i=2$,

\begin{equation*}
Q_1=(A^{(-1)}, (a_0^{(1)},A^{(1)}) ) \in K^*(1),
\end{equation*}
\begin{equation*}
Q_2=(\zeta^{-1}; B^{(-2)}, (b_0^{(-1)},B^{(-1)}), (1,\mathbf{0}) )
\in K^*(2),
\end{equation*}

\noindent where $Q_1 \ _1\infty_{-i} \ Q_2$ exists.  Using steps
(a) and (b), Propositions \ref{P:K^*_assoc} and \ref{P:H^*_assoc},
and Lemma \ref{L:doubleconvergence}

\begin{align*}
\mu_2^{\co}&(Q_1 \ _1 \infty_{-2} \ Q_2) \\
&=\mu_2^{\co}((Q_1 \ _1\infty_{-1} \ ((1,\mathbf{0}), B^{(-2)})) \
_1\infty_{-2} (\zeta^{-1}; (1,\mathbf{0}), (b_0^{(-1)},B^{(-1)}),
\mathbf{0})) \\
&=\mu_1^{\co}((Q_1) \ _1\infty_{-1} \ ((1,\mathbf{0}), B^{(-2)}))
\ _1*_{-2} \mu_2^{\co}((\zeta^{-1}; (1,\mathbf{0}),
(b_0^{(-1)},B^{(-1)}), \mathbf{0})) \\
&=((\mu_1^{\co}(Q_1) \ _1*_{-1} \ \mu_1^{\co}(((1,\mathbf{0}),
B^{(-2)}))) \ _1*_{-2} \\
& \ \ \ \ \ \ \ \ \ \ \ \ \ \ \ \ \ \ \ \ \ \ \ \ \ \ \ \ \ \ \ \
\mu_2^{\co}((\zeta^{-1}; (1,\mathbf{0}), (b_0^{(-1)},B^{(-1)}),
\mathbf{0})))
e^{-\Gamma(A^{(1)}, B^{(-2)},a_0^{(1)})d} \\
&=(\mu_1^{\co}(Q_1) \ _1*_{-2} \ (\mu_1^{\co}(((1,\mathbf{0}),
B^{(-2)})) \ _1*_{-2} \\
& \ \ \ \ \ \ \ \ \ \ \ \ \ \ \ \ \ \ \ \ \ \ \ \ \ \ \ \ \ \ \ \
\mu_2^{\co}((\zeta^{-1}; (1,\mathbf{0}), (b_0^{(-1)},B^{(-1)}),
\mathbf{0}))))
e^{-\Gamma(A^{(1)}, B^{(-2)},a_0^{(1)})d} \\
%=(\mu_1^{\co}(Q_1) \ _1*_{-2} \ \mu_2^{\co}(((1,\mathbf{0}),
%B^{(-2)}) \ _1\infty_{-2} \
%(\zeta^{-1}; (1,\mathbf{0}), (b_0^{(-1)},B^{(-1)}), \mathbf{0}))) \\
%e^{-\Gamma(A^{(1)}, B^{(-i)},a_0^{(1)})d}  \\
&=(\mu_1^{\co}(Q_1) \ _1*_{-2} \ \mu_2^{\co}(Q_2))
e^{-\Gamma(A^{(1)}, B^{(-2)},a_0^{(1)})d}
\end{align*}

Step (d): Let $i=1$,

\begin{equation*}
Q_1=(A^{(-1)}, (a_0^{(1)},A^{(1)}) ) \in K^*(1),
\end{equation*}
\begin{equation*}
Q_2=(\zeta^{-1};B^{(-2)}, (b_0^{(-1)},B^{(-1)}), (1,\mathbf{0}) )
\in K^*(2),
\end{equation*}

\noindent where $Q_1 \ _1\infty_{-i} \ Q_2$ exists.

Let $\sigma \in S_2$ be the transposition of two elements. We use
the naturality of permutations (Propositions \ref{P:K^*perm} and
\ref{P:H^*perm} and axiom 4), along with part (c) to observe the
following:

\begin{align*}
\mu_2^{\co}(Q_1 \ _1\infty_{-1} \ Q_2)
&=\mu_2^{\co}( \sigma (Q_1 \ _1\infty_{-2} \ \sigma  \ Q_2)) \\
&=\sigma (\mu_1^{\co}(Q_1) \ _1*_{-2} \ \mu_2^{\co}(\sigma  \
Q_2))
e^{-\Gamma(A^{(1)}, B^{(-1)}((b_0^{(-1)})^{-1}),a_0^{(1)})d} \\
&=\sigma (\mu_1^{\co}(Q_1) \ _1*_{-2} \ \sigma \mu_2^{\co}(Q_2))
e^{-\Gamma(A^{(1)}, B^{(-1)}((b_0^{(-1)})^{-1}),a_0^{(1)})d} \\
&=(\mu_1^{\co}(Q_1) \ _1*_{-1} \ \mu_2^{\co}(Q_2))
e^{-\Gamma(A^{(1)}, B^{(-1)}((b_0^{(-1)})^{-1}),a_0^{(1)})d}.
\end{align*}

\noindent The careful reader will have noticed the inputs for
$\Gamma$ are modified to
$A^{(1)}$, \\
$B^{(-1)}((b_0^{(-1)})^{-1})= ((b_0^{(-1)})^{-1}B^{(-1)}_1,
(b_0^{(-1)})^{-2}B^{(-1)}_2, \ldots)$, and $a_0^{(1)}$ because we
are working with $\sigma \ Q_2$ instead of $Q_2$. Via Proposition
4.2.1 in \cite{H}, we have

\begin{equation*}
\Gamma(A^{(1)}, B^{(-1)}((b_0^{(-1)})^{-1}), a_0^{(1)})
=\Gamma(A^{(1)}, B^{(-1)},a_0^{(1)}b_0^{(-1)}).
\end{equation*}

Step (e): Let $i=n$,

\begin{equation*}
Q_1=(z^{-1}; A^{(-2)}, (a_0^{(-1)},A^{(-1)}), (1,\mathbf{0}) ) \in
K^*(2),
\end{equation*}
\begin{multline*}
Q_2=(\zeta_{-n+1}^{-1}, \ldots, \zeta_{-1}^{-1};
 B^{(-n)},
(b_0^{(-n+1)},B^{(-n+1)}), \ldots, \\
(b_0^{(-1)},B^{(-1})), (b_0^{(1)},B^{(1)})  ) \in K^*(n),
\end{multline*}

\noindent where $Q_1 \ _1\infty_{-n} \ Q_2$ exists. Then the
argument mimics that of Step (b) but with
(\ref{E:firstsewingidentity}) replacing
(\ref{E:secondsewingidentity}).  (See \cite{K} for details.)

Step (f): Let $1 \leq i \leq n$,

\begin{equation*}
Q_1=(z_{-1}^{-1}; A^{(-2)}, (a_0^{(-1)},A^{(-1)}), (1,\mathbf{0})
) \in K^*(2),
\end{equation*}
\begin{multline*}
Q_2=(\zeta_{-n+1}^{-1}, \ldots, \zeta_{-1}^{-1}; B^{(-n)},
(b_0^{(-n+1)},B^{(-n+1)}), \ldots, \\
(b_0^{(-1)},B^{(-1})),  (b_0^{(1)},B^{(1)})) \in K^*(n),
\end{multline*}

\noindent such that $Q_1 \ _1\infty_{-i} \ Q_2$ exists. If
$\sigma$ is defined to be the transposition $(i \ n)$, then
following the proof of Step (d), with $\tau$ the permutation $(i \
n+1 \ i+1 \ i+2 \ldots n)$,

\begin{align*}
\mu_{n+1}^{\co}(Q_1 \ _1\infty_{-i} \ Q_2)
&=\mu_{n+1}^{\co}(\tau(Q_1 \ _1\infty_{-n} \ \sigma  Q_2)) \\
&=\tau(\mu_{n+1}^{\co}(Q_1 \ _1\infty_{-n} \ \sigma  Q_2)) \\
&=\tau(\mu_2^{\co}(Q_1) \ _1*_{-n} \ \mu_n^{\co}(\sigma \ Q_2))
e^{-\Gamma(A^{(1)}, B^{(-i)}((b_0^{(-i)})^{-1}),a_0^{(1)})d} \\
&=\tau(\mu_2^{\co}(Q_1) \ _1*_{-n} \ \sigma\mu_n^{\co}(Q_2))
e^{-\Gamma(A^{(1)}, B^{(-i)}, a_0^{(1)} b_0^{(-i)})d} \\
&=(\mu_2^{\co}(Q_1) \ _1*_{-i} \ \mu_n^{\co}(Q_2))
e^{-\Gamma(A^{(1)}, B^{(-i)}, a_0^{(1)} b_0^{(-i)})d}.
\end{align*}

\noindent where Step (e) gives the key equality.

Step (g): We will now use induction on $m$ for $m\geq 2$. Assume
that for all $\ell < m$ the sewing axiom holds and let $i \leq n$,

\begin{multline*}
Q_1=(z_{-m+1}^{-1}, \ldots, z_{-1}^{-1}; A^{(-m)},
(a_0^{(-m+1)},A^{(-m+1)}), \ldots, \\
(a_0^{(-1)},A^{(-1})), (a_0^{(1)},A^{(1)})  ) \in K^*(m),
\end{multline*}
\begin{multline*}
Q_2=(\zeta_{-n+1}^{-1}, \ldots, \zeta_{-1}^{-1}; B^{(-n)},
(b_0^{(-n+1)},B^{(-n+1)}), \ldots, \\
(b_0^{(-1)},B^{(-1})),  (b_0^{(1)},B^{(1)})) \in K^*(n),
\end{multline*}

\noindent such that $Q_1 \ _1\infty_{-i} \ Q_2$ exists.

We begin by decomposing $Q_1$ into

\begin{equation*}
Q^+_1=(z_{-m+1}^{-1}; A^{(-m)}, (a_0^{(-m+1)},A^{(-m+1})),
 (1, \mathbf{0}) ),
\end{equation*}
\begin{multline*}
Q^-_1=(z_{-m+2}^{-1}, \ldots, z_{-1}^{-1}; \mathbf{0},
(a_0^{(-m+2)},A^{(-m+2)}),  \ldots,
(a_0^{(-1)},A^{(-1})),(a_0^{(1)},A^{(1)}) ).
\end{multline*}

\noindent Then as in Step (c) (using associativity, the
$t$-contraction  and Lemma \ref{L:doubleconvergence})

\begin{align*}
\mu&_{m+n-1}^{\co} (Q_1 \ _1 \infty_{-i} Q_2) \\
&=\mu_{m+n-1}^{\co} ((Q_1^{+} \ _1\infty_{-m+1} \ Q_1^{-}) \ _1\infty_{-i} Q_2) \\
&=\mu_{m+n-1}^{\co} (Q_1^{+} \ _1\infty_{-i-m+2} \ (Q_1^{-} \ _1\infty_{-i} Q_2)) \\
&=\mu_2^{\co}(Q^+_1) \ _1*_{-i-m+2} \ \mu_{m+n-2}^{\co}
(Q^-_1 \ _1\infty_{-i} Q_2) \\
&=(\mu_2^{\co}(Q^+_1) \ _1*_{-i-m+2} \ (\mu_{m-1}^{\co}(Q^-_1) \
_1*_{-i} \mu_n^{\co}(Q_2)))
e^{-\Gamma(A^{(1)}, B^{(-i)}, a_0^{(1)} b_0^{(-i)})d}  \\
&=((\mu_2^{\co}(Q^+_1) \ _1*_{-m+1} \ \mu_{m-1}^{\co}(Q^-_1)) \
_1*_{-i} \mu_n^{\co}(Q_2))
e^{-\Gamma(A^{(1)}, B^{(-i)}, a_0^{(1)} b_0^{(-i)})d}  \\
&=(\mu_m^{\co}(Q_1) \ _1*_{-i} \ \mu_n^{\co}(Q_2))
e^{-\Gamma(A^{(1)}, B^{(-i)}, a_0^{(1)} b_0^{(-i)})d}
\end{align*}

\noindent where the third, fourth and sixth equalities employ the
inductive assumption and, in the case $m=2$, Step (f).

Step (h): We will now use induction on $n$ for $n\geq 2$. Assume
that for all $\ell < n$ the sewing axiom holds, let $i$, $Q_1$ and
$Q_2$ be as in Step (h) and decompose $Q_2$ as

%\begin{multline*} Q_1=(z_{-m+1}^{-1}, \ldots, z_{-1}^{-1}; A^{(-m)},
%(a_0^{(-m+1)},A^{(-m+1)}), \ldots, \\
%(a_0^{(-1)},A^{(-1})), (a_0^{(1)},A^{(1)})  ) \in K^*(m),
%\end{multline*}
%\begin{multline*} Q_2=(\zeta_{-n+1}^{-1} , #\ldots,
%\zeta_{-1}^{-1}; B^{(-n)},
%(b_0^{(-n+1)},B^{(-n+1)}), \ldots, \\
%(b_0^{(-1)},B^{(-1})),  (b_0^{(1)},B^{(1)})) \in K^*(n),
%\end{multline*}

%\noindent such that $Q_1 \ _1\infty_{-i} \ Q_2$ exists. This time
%we begin by decomposing $Q_2$ into

\begin{equation*}
Q^+_2=(\zeta_{-n+1}^{-1}; B^{(-n)}, (b_0^{(-n+1)},B^{(-n+1})),
(1,\mathbf{0})  ),
\end{equation*}
\begin{multline*}
Q^-_2=(\zeta_{-1}^{-n+2}, \ldots, \zeta_{-1}^{-1}; \mathbf{0},
(b_0^{(-n+2)},B^{(-n+2)}), \ldots, (b_0^{(-1)},B^{(-1})),
(b_0^{(1)},B^{(1)})  ).
\end{multline*}

\noindent There are three possibilities:

(i) If $i=n$, then

\begin{align*}
Q_1 \ _1\infty_{-n} \ Q_2 &=
Q_1 \ _1\infty_{-n} \ (Q^+_2 \ _1\infty_{-n+1} \ Q^-_2) \\
&= (Q_1 \ _1\infty_{-2} \ Q^+_2) \ _1\infty_{-n+1} \ Q^-_2.
\end{align*}

\noindent For $n>2$, the sewing follows from the inductive
assumption, associativity, the $t$-contraction  and Lemma
\ref{L:doubleconvergence}.  For $n=2$, we supplement induction
with the Steps (c) and (g).

(ii) If $i=n-1$, then

\begin{align*}
Q_1 \ _1\infty_{-n+1} \ Q_2 &=
Q_1 \ _1\infty_{-n+1} \ (Q^+_2 \ _1\infty_{-n+1} \ Q^-_2) \\
&= (Q_1 \ _1\infty_{-1} \ Q^+_2) \ _1\infty_{-n+1} \ Q^-_2
\end{align*}

\noindent For $n>2$, the sewing follows from the inductive
assumption, associativity, the $t$-contraction  and Lemma
\ref{L:doubleconvergence}.  For $n=2$, we supplement induction
with the Steps (d) and (g).

(iii) For $i<n-1$, then

\begin{align*}
Q_1 \ _1\infty_{-i} \ Q_2 &=
Q_1 \ _1\infty_{-i} \ (Q^+_2 \ _1\infty_{-n+1} \ Q^-_2) \\
&= Q^+_2 \ _1\infty_{-m-n+2} \ (Q_1 \ _1\infty_{-i} \ Q^-_2).
\end{align*}

\noindent The sewing follows from the inductive assumption, Step
(e), associativity, the $t$-contraction  and Lemma
\ref{L:doubleconvergence}.
\end{proof}

\subsection{The categorical isomorphism between GVOAs and VOAs} \label{S:iso}
Let $d$ be a complex number, $\mathbf{V}^*(d)$ be the category of
vertex operator coalgebras of rank $d$, and $\mathbf{G}^*(d)$ be
the category of geometric vertex operator coalgebras of rank $d$.
In the previous two sections we have defined functors

\begin{align*}
F_{\mathbf{V}^*(d)} :  \mathbf{V}^*(d)  \to & \mathbf{G}^*(d) \\
(V, \co, c, \rho) \mapsto & (V, \mu^{\co})&
\end{align*}
and

\begin{align*}
F_{\mathbf{G}^*(d)} :  \mathbf{G}^*(d)  \to & \mathbf{V}^*(d) \\
(V, \mu) \mapsto & (V, \co_{\mu}, c_{\mu}, \rho_{\mu}).
\end{align*}

\noindent We have shown explicitly that these two maps take
objects to objects and recalling our discussion of morphisms in
the introduction of Section \ref{S:iso1} it is evident that they
respect morphisms as well (See also \cite{K}). The main purpose of
defining $F_{\mathbf{V}^*(d)}$ and $F_{\mathbf{G}^*(d)}$ is to
show that they are inverses to each other and that the categories
of VOCs and GVOCs are isomorphic.

\begin{theorem}
The categories $\mathbf{V}^*(d)$ and $\mathbf{G}^*(d)$ are
isomorphic.  In particular, the functors $F_{\mathbf{V}^*(d)}$ and
$F_{\mathbf{G}^*(d)}$ satisfy

\begin{equation} \label{E:9231}
F_{\mathbf{G}^*(d)} F_{\mathbf{V}^*(d)}= 1_{\mathbf{V}^*(d)},
\end{equation}
\begin{equation} \label{E:9232}
F_{\mathbf{V}^*(d)} F_{\mathbf{G}^*(d)}= 1_{\mathbf{G}^*(d)}
\end{equation}

\noindent where $1_{\mathbf{V}^*(d)}$ and $1_{\mathbf{G}^*(d)}$
are the identity functors on $\mathbf{V}^*(d)$ and
$\mathbf{G}^*(d)$, respectively.
\end{theorem}

\begin{proof}
Since the second half of the theorem implies the first half, we
will simply prove (\ref{E:9231}) and (\ref{E:9232}).  First,
observe that

\begin{equation*}
c_{\mu}=\mu_0^{\co}((1,\mathbf{0}))=c,
\end{equation*}

\noindent and

\begin{equation*}
\rho_{\mu} =-\left. \frac{d}{d \epsilon} \mu_0^{\co}(1,(0,
\epsilon, 0,0, \ldots))\right|_{\epsilon=0} =-\left. \frac{d}{d
\epsilon} c e^{-\epsilon L(2)} \right|_{\epsilon=0} =cL(2) = \rho
\end{equation*}

\noindent by employing (\ref{E:cL(2)}).  For any $v' \in V \otimes
V$,

\begin{align*}
\mathrm{Res}_x x^n \langle v',\co^{\mu}(x) v \rangle
&=\mathrm{Res}_z z^n \langle v',\mu^{\co}((z^{-1} ;(1,\mathbf{0}),
(1,\mathbf{0}),\mathbf{0}))(v) \rangle \\
&=\mathrm{Res}_z z^n \left( \langle v', \co(x)v \rangle |_{x=z} \right) \\
&=\mathrm{Res}_x x^n \langle v',\co(x) v \rangle.
\end{align*}

\noindent Thus (\ref{E:9231}) is verified.

As for (\ref{E:9232}), we need only verify that

\begin{multline} \label{E:9233}
\langle v', \mu_n
(z^{-1}_{-1},\ldots,z^{-1}_{-n+1};(a_0^{(1)},A^{(1)}),
(a_0^{(-1)},A^{(-1)}),\ldots,\\
(a_0^{(-n+1)},A^{(-n+1)}),A^{(-n)})v \rangle \\
=\iota^{-1}_{x \cdots n-1} \langle v', \left( e^{-\sum_{j \in
\Z_+} A^{(-1)}_j L(-j)} (a_0^{(-1)})^{-L(0)} \otimes \cdots
\otimes e^{-\sum_{j \in \Z_+} A^{(-n+1)}_j L(-j)}
\right. \\
\left.  (a_0^{(-n+1)})^{-L(0)} \otimes e^{-\sum_{j \in \Z_+}
A^{(-n)}_j L(-j)} \right)
(\underbrace{Id_V \otimes \cdots \otimes Id_V}_{n-2} \otimes \co_{\mu}(x_{n-1})) \cdots \\
(Id_V \otimes \co_{\mu}(x_{2})) \co_{\mu}(x_{1}) \left.
e^{-\sum_{j \in \Z_+} A^{(1)}_j L(j)} ( a_0^{(1)})^{-L(0)}v
\rangle \right|_{x_{i}=z_{-i}}.
\end{multline}

Since both sides of (\ref{E:9233}) define geometric vertex
operator coalgebras, both satisfy the sewing axiom.  Therefore, by
Proposition \ref{P:gen_mod_space}, we need only prove
(\ref{E:9233}) in the specialized cases

\begin{equation*}
\langle v', \mu_2(z^{-1};
(1,\mathbf{0}),(1,\mathbf{0}),\mathbf{0})v \rangle =\langle v',
\co_{\mu}(x)v \rangle |_{x=z},
\end{equation*}
\begin{equation*}
\mu_0((1,\mathbf{0}))v = c,
\end{equation*}
\begin{equation*}
\langle v', \mu_1((a_0^{(1)},A^{(1)}), A^{(-1)})v \rangle =\langle
v', e^{-\sum_{j \in \Z_+} A^{(-1)}_j L(-j)} e^{-\sum_{j \in \Z_+}
A^{(1)}_j L(j)} ( a_0^{(1)})^{-L(0)}v \rangle.
\end{equation*}

The first 2 of these equations are true by the definition of
$F_{\mathbf{G}^*(d)}$. The third equation, the proof for $K^*(1)$,
is identical to the proof for $K(1)$ in \cite{H} (Equation
(5.4.31)) since $K(1)=K^*(1)$ in the most specific sense, the
definitions of $L(0)$ are the same and the meromorphicity axiom is
the same in that special case.

%%%%%If I wanted to write out the proof of the third equation, this is the start.
%For the third equation, we can again view it as the sewing of
%\begin{equation*}
%\langle v', \mu_1((1,\mathbf{0}), A^{(-1)})v \rangle
%=\langle v', e^{-\sum_{j \in \Z_+} A^{(-1)}_j L(-j)} v \rangle.
%\end{equation*}
%\begin{equation*}
%\langle v', \mu_1((1,A^{(1)}), \mathbf{0})v \rangle
%=\langle v', e^{-\sum_{j \in \Z_+} A^{(1)}_j L(j)} v \rangle.
%\end{equation*}
%\begin{equation*}
%\langle v', \mu_1((a_0^{(1)},\mathbf{0}), \mathbf{0})v \rangle
%=\langle v', ( a_0^{(1)})^{-L(0)}v \rangle.
%\end{equation*}

\end{proof}

\end{document}